\documentclass[]{article}

\usepackage{hyperref} 
\usepackage{amssymb,amsmath}
\usepackage{graphicx}
\usepackage{color}
\usepackage{overpic}
\usepackage{tikz}

\newcommand\blfootnote[1]{%
  \begingroup
  \renewcommand\thefootnote{}\footnote{#1}%
  \addtocounter{footnote}{0}%
  \endgroup
}

\title{\bf 
Dissolution of Variable-in-Shape Drug Particles via the Level-Set Method
}
\author{
\phantom{}\blfootnote{Istituto per le Applicazioni del Calcolo, Consiglio Nazionale delle Ricerche, Rome, Italy.}
\phantom{}\blfootnote{Dipartimento di Ingegneria e Architettura, Universit\`a di Trieste, Trieste, Italy}
\phantom{}\blfootnote{Dipartimento di Scienze di Base e Applicate per l'Ingegneria, Sapienza Universit\`a di Roma, Rome, Italy}
Emiliano Cristiani$^*$, 
Mario Grassi$^\dag$, \\ 
Francesca L. Ignoto$^{*,\ddag}$, 
Giuseppe Pontrelli$^*$
}
\date{\today}

\newtheorem{rem}{Remark}

\renewcommand{\u}{\varphi}
\newcommand{\R}{\mathbb R}
\newcommand{\Gv}{\mathbf v}
\newcommand{\Gx}{\mathbf x}
\newcommand{\Gq}{\mathbf q}
\newcommand{\Gn}{\mathbf {\hat n}}
\newcommand{\Grad}{\mathbf \nabla}
\newcommand{\rosc}{R}

\newcommand{\req}{R_{\textsc{eq}}}
\newcommand{\rplane}{R_{\textsc{plane}}}
\newcommand{\V}{V_{\textsc{ext}}}
\newcommand{\be}{\begin{equation}}
\newcommand{\ee}{\end{equation}}
\newcommand{\ds}{\displaystyle}

\begin{document}

\maketitle

\begin{abstract} 
In this work, we deal with a mathematical model  describing the dissolution process of irregularly shaped particles. 
In particular, we consider a complete dissolution model accounting for both surface kinetics and convective diffusion, for three drugs with different solubility and wettability: theophylline, griseofulvin, and nimesulide.  
The possible subsequent recrystallization process in the bulk fluid is also considered. 
The governing differential equations are revisited in the context of the level-set method and Hamilton-Jacobi equations, then they are solved numerically.
This choice allows us to deal with the simultaneous dissolution of hundreds of different polydisperse particles.
We show the results of many computer simulations which investigate the impact of the particle size, shape, area/volume ratio, and the dependence of the interfacial mass transport coefficient on the surface curvature.
%
\end{abstract}


\section{Introduction}\label{sec:intro}
The dissolution of a solid (drug) in a liquid phase can be defined as ``the mixing of two phases (the solid and the liquid one) with the formation of a new homogeneous phase (i.e., the solution)'' \cite{siepmann2013}. 
Thus, dissolution implies the destroying of the original crystalline network and the subsequent delivery of ``solid'' atoms/molecules in the liquid phase. Kinetically speaking, it can be considered as the combination of five concurrent steps \cite{siepmann2013,hasa2014,parmar2018}, 
comprising: 
\begin{enumerate}
    \item the contact of the liquid with the solid surface (wetting), giving origin to the formation of a new solid-liquid interface starting from the solid-vapor one; 
    \item the breakdown of the crystalline network (fusion);
    \item the transfer of solid molecules to the solid-liquid interface (solvation); 
    \item the diffusion of the solvated molecules through the unstirred boundary layer surrounding the solid surface (diffusion);
    \item the convective transport of solvated drug molecules into the well stirred bulk solution (convection).
\end{enumerate}
As each one of the first four steps implies to overcome an energetic barrier, the sum of these energies represents the total resistance to solid (drug) dissolution. Obviously, the higher the dissolution energy required (i.e., the higher the mass transfer resistance) the lower the dissolution kinetics is. 

The analysis of the dissolution phenomenon is not only scientifically interesting \textit{per se}, but it is absolutely relevant for the pharmaceutical field. In fact, it affects the bioavailability of drugs, i.e., the rate and extent to which the drug is absorbed from a delivery system and becomes available at the site of drug action \cite{kwan1997,naidu2008}. 
Bioavailability depends on both the drug permeability through the cell membrane and on the drug dissolution in physiological fluids, as solid drugs cannot be absorbed by living tissues \cite{kwan1997}. 
This last aspect is very important remembering that about $40\%$ of the market drugs and $70-90\%$ of new chemical entities are characterized by low dissolution kinetics due to their poor water solubility \cite{loftsson2010,davis2018,gigliobianco2018,bertoni2019}.
Accordingly, scientific and industrial reasons have attracted the attention of many researchers on dissolution since very long times. 
The seminal work of Hixson and Crowell \cite{hixson1931a,hixson1931b,hixson1931c}, who, for the first time, considered the effect of surface reduction upon the dissolution of spherical particles and established the well-known cubic law, dates back to 1931. 
Later on, Niebergall et al.\ \cite{niebergall1963} noted a deviation from the cubic law and proposed an improvement assuming that the thickness of the diffusion layer surrounding the dissolving particle was proportional to the square root of the mean particle diameter. 
The elegant approach of Pedersen and coworkers extended the mathematical modelling of dissolution to polydisperse spherical particles \cite{pedersen1975a, pedersen1975b, pedersen1977, pedersen1978}. 
Notably, this model reduces to the Hixson-Crowell model in the case of monodisperse systems. 
Other authors argued that the overall dissolution process can be affected by the occurrence of a surface reaction between solute and solvent molecules or by limited solid surface solubility \cite{grassi2007}. 
In either case, the final result is a time-dependent drug concentration at the solid-liquid interface that is lower than the drug solubility in the solvent. 
A similar phenomenon concerns metastable solids undergoing a phase change (amorphous-crystalline or polymorphic transformation) upon dissolution, this reflecting in a time-dependent solid drug solubility \cite{guo2018,abrami2020ADMET}.
Interestingly, the possible drug degradation in the bulk fluid after dissolution was also considered \cite{thormann2014}. 

Another aspect that attracted the interest of researchers was the effect of particles shape on dissolution \cite{mosharraf1995}. 
While Hirai et al.\ \cite{hirai2017} did not explicitly consider the shape of the particles but focused instead on a law which can describe the time dependence of the dissolution surface, Abrami et al.\ \cite{abrami2020ADMET} dealt with spherical, cylindrical and parallelepiped particles, Hsu and Wu \cite{hsu2021} considered spheres, cylinders, bi-cones, cones and inverse-cones shaped particles, and Yuan et al.\ \cite{yuan2013} focused on the dissolution of irregularly shaped particles. 
It is easy to predict that the symbiotic-synergic effect occurring between these (and other) theoretical models and the emerging experimental set up aimed to visualize particle shape variation upon dissolution \cite{abrami2020ADMET,lausch2024} will soon lead to important new achievements in this field.

Undoubtedly, particle shape is also connected with two other important aspects, i.e.\ drug concentration profile in the boundary layer (BL) surrounding the solid surface, and BL thickness. 
Indeed, drug concentration profile is not linear, as originally assumed \cite{hsu1991,wang1999,hasa2014,salish2024}, unless dissolution occurs from a flat surface. 
In addition, BL thickness depends on both dissolution medium hydrodynamic conditions and on particle dimension, as nicely documented by D’Arcy and Persoons \cite{darcy2011,darcy2019}.

To the best of our knowledge, one of the most recent and comprehensive mathematical model accounting for many aspects connected to dissolution and presented in this introduction is that by Abrami et al.\ \cite{abrami2022EJPB} who employed the concept of local curvature radius to describe the dissolution of concave and convex surfaces of any shape. 
The target of this work is to move forward, developing a mathematical model capable to account also for the simultaneous presence of solid particles with different shapes and dimensions, this being one of the most common situation met in real conditions. 
Moreover, the effect of particles surface/volume ratio (or perimeter/surface ratio in the two dimensional case), and the curvature radius dependence of the interfacial mass transport coefficient will be considered.
To this purpose, our simulations dealt with three widely used drugs, namely \emph{theophylline} (soluble and wettable-bronco-dilatator drug), \emph{griseofulvin} (poorly soluble but wettable-antifungal drug) and \emph{nimesulide} (poorly soluble and poorly wettable-non-steroidal anti-inflammatory drug).

The mathematical model of the dissolution dynamics is revisited in the context of the level-set method and Hamilton--Jacobi equations \cite{osher1988JCP}, which revealed very useful and versatile as it allows to manage hundreds of particles with a reasonable computational effort (computing solution requires some minutes). 
Moreover, it allows to easily manage the presence of ``corners'' in the particles shape. 
All the simulations were performed in two dimensions but the extension to three dimensions is also possible.

\section{Modelling drug dissolution}\label{sec:modelling}
As discussed in Sect.\ \ref{sec:intro}, dissolution is a complex process made up by a series of different sub-phenomena, each one representing an energetic barrier. Accordingly, the inverse of the overall dissolution resistance equals the sum of the resistance inverse competing to each sub-phenomenon \cite{siepmann2013}. 
The aim of this section is to translate into mathematical terms the above mentioned sub-phenomena in order to get a model able to describe the time variation of drug concentration inside the release environment.

In this paper we overcome the usual assumption of particles of spherical shape in the attempt to understand the effect of the local curvature on the dissolution process. 
Since particles change their shape and size during dissolution, this reflects in a change of dissolution rate over time. 

Let us consider a particle $P_0$ of arbitrary shape in $\R^3$ with initial volume $V_0$, enclosed by a surface $S_0$. 
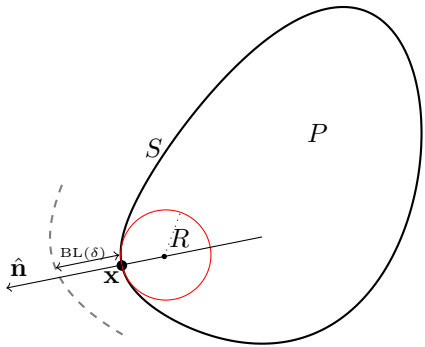
\begin{figure}[t!]
\centering
\begin{tikzpicture}[scale=2]
	\draw[black,thick,domain=0:360,samples=100,smooth] 
	plot(
	{
		(1 + 0.06 * cos(\x)) * cos(\x)
	},
	{
		(1 + 0.4 * cos(\x)) * (sin(\x)+1) - 0.4
	});
	
	\draw[dashed,gray,thick,domain=160:230,samples=100,smooth] 
	plot(
	{
		1.5*((1 + 0.06 * cos(\x)) * cos(\x))
	},
	{
		1.5*((1 + 0.4 * cos(\x)) * (sin(\x)+1) - 0.4)
	});
	
	\fill (-0.65,0.18) circle (0.5pt); 
	\fill (-0.932,0.12) circle (1pt); 
	\draw[red] (-0.64,0.19) circle (0.3); 
	\draw[dotted] (-0.65,0.18) -- (-0.54,0.48); 
	
        \draw[<-] (-1.7,-0.03) -- (0,0.31); 
	
	\draw[<->] (-1.37,0.105) -- (-0.95,0.19);
	
	\node[] at (-1,0.04) {$\mathbf x$};
	\node[] at (-1.19,0.2) {\tiny BL($\delta$)};
	\node[] at (-0.55,0.3) {$R$};
        \node[] at (-0.72,0.9) {$S$};
	\node[left] at (0.5,1) {$P$};
	\node[left] at (-1.5,0.11) {$\hat{\mathbf n}$};
	
\end{tikzpicture}

	
	
	
%
	
\caption{
Two-dimensional representation of a particle $P$ dissolving in the bulk fluid, with osculating circle of radius $R$ at a generic point $\Gx$ on its surface $S$. Local dissolution occurs in normal direction with respect to $S$. 
}
\label{fig:basics}
\end{figure}
The particles dissolves in an external bulk fluid of constant volume $\V \gg V_0$. In the following, we will denote by $V^+=\V / V_0$ the ratio between the two volumes.
Moreover, an interaction between the solute and solvent (solvation) at solid-liquid interface yields the formation of a BL of thickness $\delta$ surrounding the solid surface \cite{abrami2022EJPB}.
During dissolution, the surface evolves and volume reduces accordingly.
In the following we will denote by $P(t)$ the dissolving particle at time $t$, by $V(t)$ its volume, and by $S(t)$ its external surface, see Fig.\  \ref{fig:basics}.

To evaluate the amount of drug dissolved, the mass balance principle is used.
By considering an infinitesimal part of the surface $dS(\Gx)$ in a neighborhood of the point $\Gx \in S$,\footnote{For the sake of simplicity, in the following we omit the dependence of the variables on $\Gx$ and $t$, unless required by the context.}
the corresponding loss of mass 
equals the drug flux leaving $dS$ \cite{abrami2022EJPB}: 
\be\label{eq31}
  \frac{dM^b}{dt} = -\frac{dM^s}{dt} = \V \frac{dC_b}{dt} = \int_S
   K(\Gx) (C_s-C_b) dS ,
\ee
 where $\rho_s$ is the particle density, $C_b(t)$ is the concentration of the bulk fluid, and $K(\Gx)$ is the position-dependent mass transfer coefficient (see dedicated Sect.\ \ref{sec:K} later on) \cite{abrami2020ADMET, abrami2022EJPB}.
Eq.\ \eqref{eq31} states that in the neighborhood of $\Gx$ the solid mass reduces of $K(\Gx) (C_s-C_b) dS$  in the time unit and,
correspondingly, the same amount in transferred into the fluid. 
In other words, \eqref{eq31} shows that the dissolution rate is position-dependent, with the implicit assumption that dissolution at $\Gx$ occurs in the \emph{normal direction} with respect to the surface. 

In \eqref{eq31} the solubility $C_s(t)$ is an exponentially decreasing function of time, as a consequence of a possible drug phase transformation at the solid-liquid interface \cite{nogami1969}:
\be\label{eq50}
C_s(t)=C_{sf}+(C_{s0}-C_{sf})\exp(-k_r t),
\ee
where $C_{s0}$ and $C_{sf}$  ($C_{s0}\geq C_{sf}$) are the initial and final values of drug solubility, respectively, and $k_r$ is the surface recrystallization constant. 
Should $C_b(t)$ exceeds $C_s(t)$, drug dissolution ends and drug re-crystallization in the release environment starts. 
This phenomenon can be properly described by the well-known Nogami model \cite{nogami1969}: 
\be \label{eq37} 
\frac{dM_c}{dt}=-k_{rb} \V (C_s - C_b),
\ee
with $M_c$ the recrystallized mass and $k_{rb}$ the bulk recrystallization constant \cite{nogami1969}. 
At every step of the process, the conservation of mass holds, that is:
\be \label{eq12}
M_0= M_s(t)+C_b(t) \V +M_c(t),
\ee
with $M_0$ the initial particle mass, $M_s$ the particle mass and $C_b(t) \V$ the drug mass dissolved in the bulk fluid.
By time differentiating \eqref{eq12}, we have
\be
\frac{d}{dt}\Big(M_s+C_b \V + M_c \Big)=0,
\ee
 and summing up the contributions given by \eqref{eq31},
\be
\V \frac{dC_b}{dt} + \frac{dM_c}{dt}= (C_s-C_b) \int_S  K(\Gx) dS  
\ee
i.e. 
\be
\frac{dC_b}{dt} = \frac{(C_s-C_b)\ds \int_S  K(\Gx) dS  -  \ds \frac{dM_c}{dt} }{ \V } . \label{eq38}
\ee
By virtue of \eqref{eq37}, \eqref{eq38} is finally replaced by 
\be
\frac{d C_b}{dt} = \frac{C_s-C_b}{\V}\int_S  K(\Gx) dS  + k_{rb}   (C_s-C_b).  
\label{eq34} 
\ee
The last equation says that the variation of concentration $C_b$ is 
due in part to the dissolution of the solid (positive contribution) and in part  to  the recrystallization (negative contribution). Actually, dissolution and recrystallization are not simultaneous, and
the whole process can be subdivided in two subsequent regimes:

\paragraph{1) Dissolution.} 
In the first phase, as long as $C_b < C_s$ the dissolution only takes place, and $C_b$ increases. 
The dissolution process is described
by the following equations:
\begin{equation}\label{eq29}
\left\{
\begin{array}{l}
\displaystyle \frac{d\Gx}{d t}  =  - \frac{K(\Gx)}{\rho_s} (C_s-C_b)\hat{\mathbf n}  \\ [4mm]
\displaystyle \frac{d C_b}{dt} = \frac{C_s-C_b}{\V} \int_S  K(\Gx) dS
\end{array}
\right., 
\qquad C_b < C_s ,
\end{equation}
with $\hat {\mathbf n}$ the unit external normal (note the ODEs are nonlinear because of $K(\Gx)$),
and initial conditions
\be
\label{initial}
\Gx(0)=\Gx_0,   \qquad C_b(0)=0,      
\ee
for any $\Gx_0 \in S_0$. 

\begin{rem}[spherical particles]\label{rem:cerchi}
If $P_0$ is a sphere, the particle maintains its original shape during dissolution and the computation simplifies. 
Therefore, the dynamics of $P(t)$ can be described by the evolution of any point on its surface or, equivalently, of its radius $R$. 
Moreover, $K$ is constant  and \eqref{eq29}-\eqref{initial} simplifies  as:
\begin{equation}\label{mainequation_singlecircularparticle}
\left\{
\begin{array}{l}
\displaystyle\frac{d R }{d t}=-\displaystyle\frac{K}{\rho_s} \big(C_s-C_b\big)\hat {\mathbf n}  \\ [4mm]
\displaystyle\frac{dC_b}{dt}=4\pi\rosc^2 K \, \frac{C_s-C_b}{\V} 
\end{array}
\right. ,
\qquad C_b<C_s,
\end{equation}
\end{rem}
with
\be
R(0)=R_0,   \qquad C_b(0)=0.      
\ee

\paragraph{2) Recrystallization.} 
In the second step, which starts at the time $t^*$ such that $C_b(t^*) = C_s(t^*)$, dissolution stops and, as long as $C_b \geq C_s$, the recrystallization phase takes place with a progressive reduction of $C_b$. This lasts until the time $t_f$ such that the equivalence $C_b(t_f)=C_s(t_f)$ holds true for the second time. 
This process is governed by the following linear equations:
\begin{equation}\label{eq21}
\left\{
\begin{array}{l}
\displaystyle \frac{dM_c}{dt} = -k_{rb} \V (C_s - C_b) \\ [4mm]
\displaystyle \frac{d \Gx}{d t} = 0 \\ [4mm]
\displaystyle \frac{d C_b}{dt}=  k_{rb} (C_s-C_b)   \end{array}
\right., 
\qquad C_b \geq C_s ,
\end{equation}
with initial conditions:
\be
M_c(t^*)=0,   \qquad 
C_b(t^*)=C_s(t^*).      \label{eq43}
\ee
Eqs.\ \eqref{eq21}-\eqref{eq43} admit the following analytical solution:
\begin{multline}
C_b(t)=C_s(t^*)\exp(-k_{rb} (t-t^*)) + \\ + (C_{s0}-C_{sf})\frac{k_{rb}}{k_{rb}-k_r}\exp(-k_{rb}t)\Big(\exp((k_{rb}-k_r)t)-\exp((k_{rb}-k_r)t^*)\Big) + \\ + C_{sf}(1-\exp(-k_{rb}(t-t^*))) 
\end{multline}
%
%
and
\begin{multline}
\frac{M_c(t)}{\V \, k_{rb}}= 
-\frac{C_s(t^*)}{k_{rb}}\exp(k_{rb}t^*)\Big( \exp(-k_{rb}t) - \exp(-k_{rb}t^*) \Big) + \\ +
\frac{C_{s0}-C_{sf}}{k_r}\left( 1- \frac{k_{rb}}{k_{rb}-k_r} \right) 
\Big( \exp(-k_r t) - \exp(-k_r t^*) \Big) + \\ +
\frac{k_{rb}(C_{s0}-C_{sf})\exp((k_{rb}-k_r)t^*)+(k_{rb}-k_r)C_{sf}\exp(k_{rb}t^*)}{k_{rb}(k_{rb}-k_r)} \Big( \exp(-k_{rb}t) - \exp(-k_{rb}t^*) \Big).
\end{multline}


%
%
%
%
%
%
%

\section{The mass transfer coefficient $K$}\label{sec:K}
As discussed in \cite{abrami2022EJPB}, dissolution depends on solid wettability with respect to the dissolving fluid, on the surface curvature, and on the thickness of the BL, which -- in turn -- depends on the surface curvature and on the relative velocity between solid surface and fluid \cite{darcy2011}.

Assuming that 1) dissolution occurs in the normal direction with respect to the solid surface, 2) mass transport inside the BL occurs according to a one-dimensional diffusive mechanism, and 3) pseudo-stationary conditions are rapidly reached in the BL, it can be demonstrated that overall mass transfer coefficient $K$ is given by
\begin{equation}\label{K}
   K=\frac{\ds\frac{1}{\sigma}}{\displaystyle\frac{1}{k_d}+\displaystyle\frac{\sigma}{k_m}}= \left[ \sigma \left(\displaystyle\frac{1}{k_d}+\displaystyle\frac{\sigma}{k_m}\right) \right]^{-1}
\end{equation}
where $k_d$ is the hydrodynamic mass transfer coefficient, connected to the BL local thickness, $k_m$ is the interfacial mass transfer coefficient, and $\sigma$ denotes the ratio between the maximum and minimum bounds of the BL. 
For further details the reader is referred to \cite{abrami2020ADMET}.

\paragraph{Curvature radius.}
It is well known that at any point $\Gx$ of a surface $S$ one can define two principal curvatures $\kappa_1$ and $\kappa_2$, defined as the maximum and minimum of the curvature of a curve contained in the surface and passing through the point $\Gx$. With these curvatures, one can then define the \textit{mean curvature} $\kappa_M=(\kappa_1+\kappa_2)/2$.
For our purposes, we define a suitable indicator of the curvature radius at any point $\Gx$ of the surface as:
$$
R=\frac{1}{\kappa_M}. 
$$
Note that, in 2D, the curvature $\kappa$ is uniquely defined, and  $R= \ds\frac{1}{\kappa}$ corresponds to the radius of the osculating circle, see Fig.\ \ref{fig:basics}. As a consequence, $\sigma$ in \eqref{K} can be written as a function of $R$, as
\be\label{def:sigma}
\sigma = \ds\frac{R+\delta}{R}=1+\frac{\delta}{R}=1+\frac{D}{k_d \, R},
\ee
where we have used the known relation $k_d = D/\delta$ \cite{levich1962book}.

\medskip

In previous papers it has been shown, by experimental evidence, that the hydrodynamic mass transfer (or intrinsic dissolution) $k_d$ is inversely related to the radius for spherical particles \cite{darcy2011}. 
The local dependence of $k_d$ is extended to any surface, through the following expression \cite{abrami2022EJPB}: 
\be
\label{kd}
   k_d= \ds\frac{D}{2R}\left(2+0.6\sqrt{\frac{2R \, \Delta U }{\nu_f}}\left(\frac{\nu_f}{D}\right)^{1/3}\right),
\ee
where $D$ is the drug diffusivity in the bulk fluid, $\nu_f$ is the kinematic viscosity, and 
\begin{equation}\label{DeltaU}
   \Delta U=(\rho_s-\rho_f)g \displaystyle\frac{(2\req)^2}{18\eta_f}
\end{equation}
is the relative solid-liquid velocity (with $\rho_f$ the fluid density, $\eta_f$ the fluid dynamic viscosity, $g$ the gravity acceleration) \cite{abrami2020ADMET}. 
In \eqref{DeltaU}, the time-dependent quantity $\req$ is defined as the radius of an {\em equivalent sphere} having the current particle volume
\begin{equation}\label{req}
   \req=\sqrt[3]{\frac{3V}{4\pi}}.
\end{equation}

Note that, in the special case of a spherical particle, $\req(t)=R(t)$ at any time $t$. 
Note also, in \eqref{kd}-\eqref{DeltaU}, the nonlinear dependence of $k_d$  on the particle size via $\req$. For an almost flat surface ($R \rightarrow \infty$), (\ref{kd}) ceases to apply and we resort to the Levich approach \cite{levich1962book}:
\be
\label{kd_flat}
 k_d=0.621 \ D^{2/3} \nu_f^{-1/6}\sqrt{\ds\frac{\Delta U}{\req}}.
\ee

As \eqref{kd} and \eqref{kd_flat} provide approximately the same $k_d$ evaluation for sufficiently wide curvature radii, following the strategy proposed in \cite{abrami2022EJPB} we define a limiting radius $\rplane$ beyond which \eqref{kd_flat} replaces \eqref{kd}: 
$\rplane$ is simply determined by equating the two equations and then solving for the unknown curvature radius, whose analytical solution reads:
\begin{align}\label{xiplane}
& \rplane=\frac{\gamma_1}{2}\left( \frac{\gamma_1+\sqrt{\gamma_1^2+4D\gamma_2}}{\gamma_2^2} \right)+\frac{D}{\gamma_2}, \\
&   \gamma_1=0.1 D^{2/3}\nu_f^{-1/6}\sqrt{\frac{4g(\rho_s-\rho_f)}{\eta_f}}\req \nonumber, \\
&   \gamma_2=0.207 D^{2/3}\nu_f^{-1/6}\sqrt{\frac{2g(\rho_s-\rho_f)}{\eta_f}}\sqrt{\req}. \nonumber
\end{align}
Eq.\ \eqref{xiplane} represents a simple strategy to manage the transition from a convex ($R\leq\rplane$) to a planar ($R>\rplane$) dissolving surface.
It is worth noting that $\rplane$ depends on the particle size via $\req$.

\medskip

Regarding instead the interfacial mass transfer coefficient $k_m$, it depends on the local wetting properties of the dissolution surface. 
In this work we overcome the common hypothesis 
of a constant $k_m$, and adopt the following curvature-dependent expression:
\begin{equation}\label{km}
k_m=k_m^\infty\left(\displaystyle\frac{\alpha}{R^3}+\displaystyle\frac{R}{R+2d_T}\right),
\end{equation}
with $k_m^\infty$ the interfacial mass transfer coefficient associated with a planar surface ($R \rightarrow \infty$), $d_T$ the Tolman length and $\alpha$ a fitting parameter. 
Likewise $k_d$, $k_m$ is also position-dependent via $\rosc$. See Appendix \ref{app:derivation_km} for the theoretical derivation of \eqref{km}.

\medskip 

\begin{figure}[h!]
    \centering
    \includegraphics[width=0.49\textwidth]{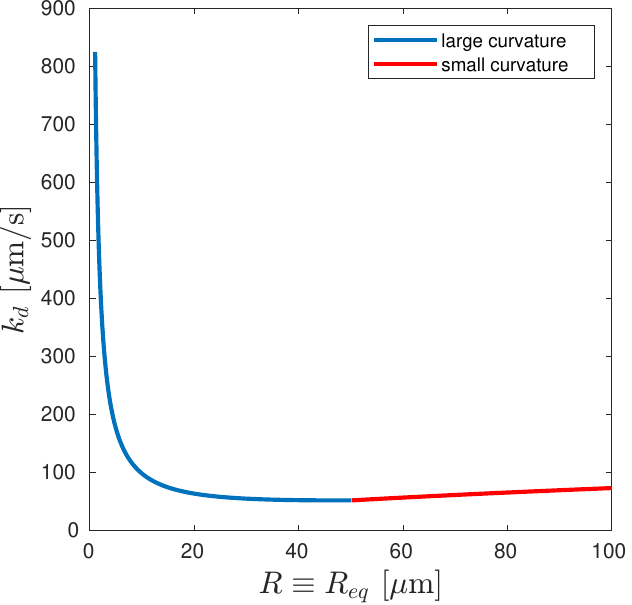}
    \includegraphics[width=0.49\textwidth]{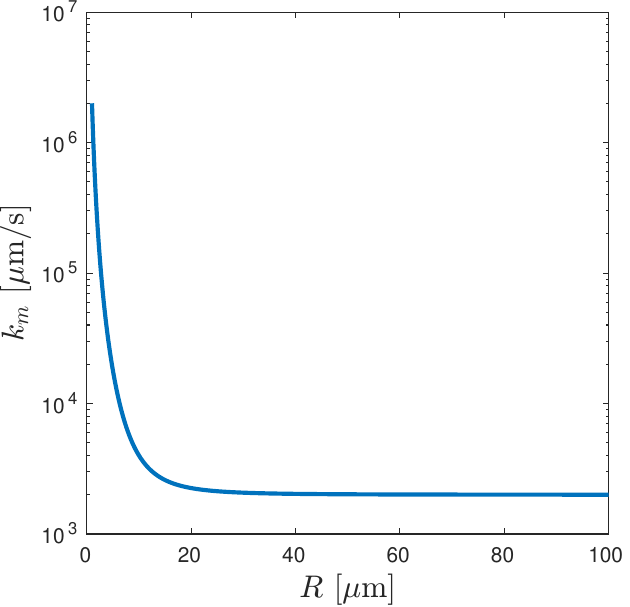}
    \caption{
    Behaviour of $ k_d(\rosc)$ (left) and $k_m(\rosc)$ (right) for a spherical particle of theophylline 37$^\circ$C of initial radius 100 $\mu$m shrinking to 0 upon complete dissolution (other parameters in Table \ref{tab:parameters}). 
    Different colors highlights the switch between the two expressions of $k_d$ \eqref{kd} and \eqref{kd_flat}.
    }
    \label{fig:funzioni_analitiche_teofillina37}
\end{figure}

Fig.\ \ref{fig:funzioni_analitiche_teofillina37} shows the graphs of $k_d$ and $k_m$ as a function of $\rosc$ for a fully-dissolving spherical particle of initial radius $R_0=100$ $\mu$m. 
It is clear that both mass transfer coefficients increase for smaller $R$. However, while $k_m$ is monotonic decreasing with $R$, $k_d$ shows a minimum. 

\begin{figure}[h!]
    \centering
    \includegraphics[width=0.49\textwidth]{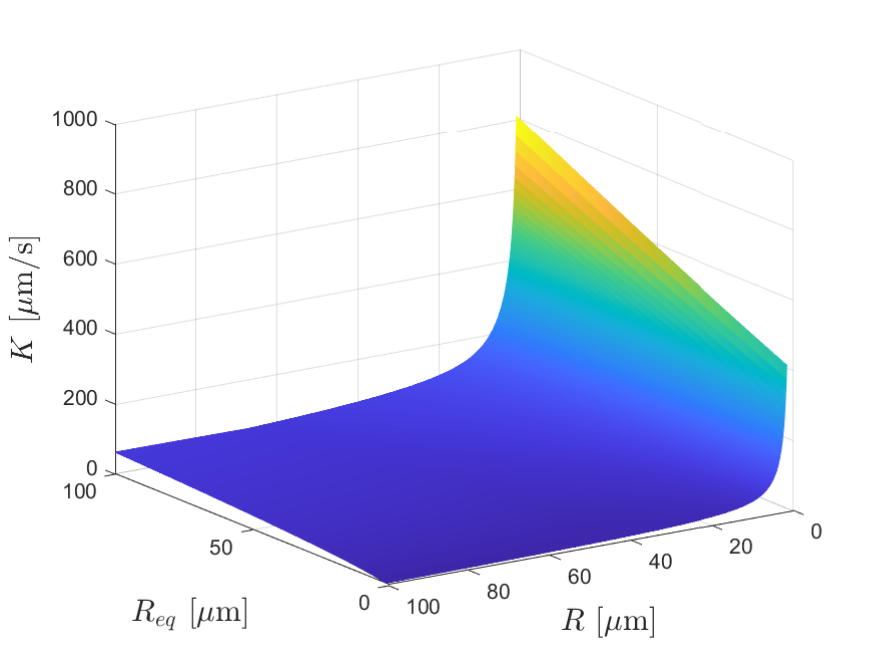}
    \includegraphics[width=0.49\textwidth]{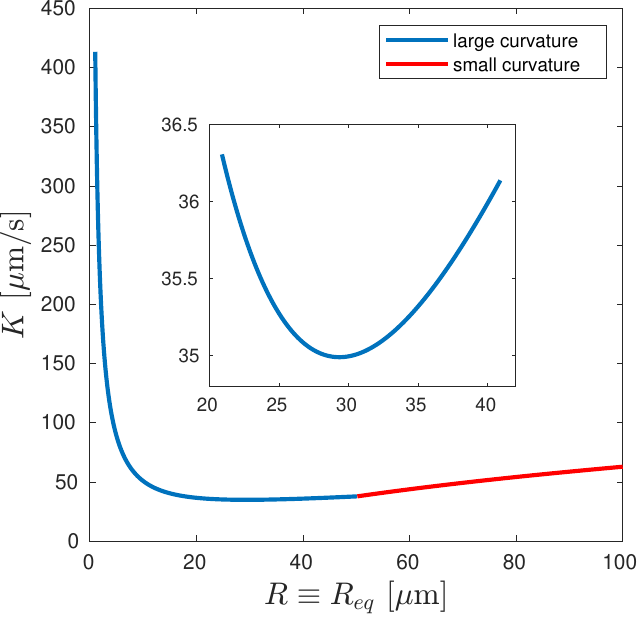}
    \caption{
    Behaviour of $K(\rosc,\req)$ (left) and $K(\rosc,\req=\rosc)$ (right, for a sphere) for a particle of theophylline 37$^\circ$C (other parameters in Table \ref{tab:parameters}). Note that $K(\rosc,\rosc)$ is not monotone (see inset graph).
    Different colors highlights the switch between the two expressions of $k_d$ \eqref{kd} and \eqref{kd_flat}.
    }
    \label{fig:funzioni_analitiche_teofillina37_K}
\end{figure}

Fig.\ \ref{fig:funzioni_analitiche_teofillina37_K} shows the graphs of $K$ as a function of $\rosc$ and $\req$ either in the general case and in the particular case of a spherical particle of initial radius $R_0=100$ $\mu$m.
It is clear that the dependence of $K$ on $R$ is more similar to that of $k_d$ rather than to that of $k_m$. 
This is simply due to the fact that the drug considered for these plots (theophylline) is a water wettable substance so that $k_m$ is much higher than $k_d$. 
Obviously, for poorly water wettable drugs, the predominance of the $k_m$ on $K$ takes place. 
All these considerations make clear the complex relation occurring between $K$ and $R$ also for spherical particles where $R=\req$.

\clearpage 
\section{The level-set methodology}\label{sec:LS}

The level-set method was introduced in \cite{osher1988JCP} and since then it was successfully applied in many contexts, see e.g., \cite{osherbook, sethianbook}. 
It allows to track Eulerianly the evolution of a $(d-1)$-dimensional surface $S$ embedded in $\R^d$ and transported by a given velocity vector field $\Gv(t,\Gx;S)$ possibly dependent on time, space, and also on the shape of the surface itself. 
In our specific case, $S$ is the boundary of the drug particle, and $\Gv$ controls its motion.  
Using the level-set method, every point $\Gx$ of the surface is moved according to the velocity field $\Gv(t,\Gx)$. 
Note that, using this approach, is \emph{not} possible tracking Lagrangianly every single point of the particle surface from the beginning to the end of its evolution, as it is done in \cite{abrami2022EJPB}. 

Let us briefly recall the method in the case $d=3$.
The main idea of the level-set method stems on the definition of a level-set function 
$\u(\Gx,t):\R^3 \times \R^+\to\R$ such that
\begin{equation}\label{Sigmat}
S(t)=\{\Gx: \u(\Gx, t)=0\},\quad \forall \, t\geq 0.
\end{equation}
In this way the surface is recovered as the zero level set of $\u$ at any time. 
Initially, the function $\u$ is chosen in such a way that
\begin{equation}\label{phi0}
\u(\Gx,0)=\u_0(\Gx) 
\left\{
\begin{array}{ll}
>0, & \textrm{if } \Gx\notin \overline{P_0}, \\
=0, & \textrm{if } \Gx\in S_0, \\
<0, & \textrm{if } \Gx\in P_0.
\end{array}
\right.
\end{equation}
A typical choice for $\u(\Gx,0)$ is the signed distance function from $S_0$, although this choice does not lead to a smooth function. 
One can prove \cite{sethianbook} that the level-set function $\u$ satisfies the following time-dependent Hamilton--Jacobi equation
\begin{equation}\label{advec_eq}
\frac{\partial \u}{\partial t}+\Gv\cdot \Grad\u=0,\qquad \Gx\in\R^3, \quad  t\in\R^+,
\end{equation}
with a suitable initial condition $\u(\Gx,0)=\u_0(\Gx)$ satisfying (\ref{phi0}). 
If the vector field $\Gv$ always points in the direction normal to the surface, i.e.\ it has the form $\Gv = v\Gn$, where $\Gn$ is the unit exterior normal and $v$ is some scalar function, the equation (\ref{advec_eq}) turns into the time-dependent Eikonal equation
\begin{equation}\label{LSequation}
\frac{\partial \u}{\partial t}+v|\Grad\u|=0,\qquad \Gx\in\R^3, \quad  t\in\R^+.
\end{equation}
This is indeed our case, since the the drug particles dissolve in normal direction.

\medskip 

One of the most appealing features of the level-set method is that several geometrical properties of the evolving surface can be described by means of its level-set function $\u$. 
For example, it possible to express the unit exterior normal $\Gn$, the mean curvature $\kappa_M$, $ R, \, S$ and $V$ in terms of $\u$ and its derivatives \cite{sethianbook}. More precisely, we have
\begin{equation}\label{levelsetvariables}
\Gn=\frac{\Grad\u}{|\Grad\u|}, \quad 
\kappa_M= \text{div } \Gn, \quad 
\rosc=\frac{1}{\kappa_M},\quad 
S=\int_{\{\phi=0\}}ds, \quad
V=\int_{\{\phi<0\}}d\Gx.
\end{equation}
\begin{rem}
    Note that all the quantities introduced so far must be are well defined everywhere in the space $\R^3$ and not only at the surface points (although they could have a clear physical meaning only there). 
    This is a crucial fact for using the level-set method, since the evolution of the surface is actually replaced by the evolution of the level-set function $\u$, which is defined on all over the space and represents a sort of extension of the surface.
\end{rem}

\subsubsection*{Tracking the particle dissolution}

Reformulating \eqref{eq29} in the framework of the level-set method, we get that the evolution of the surface is described by the following dynamics:
\begin{equation}\label{mainequation_singleparticle}
\left\{
\begin{array}{l}
\displaystyle\frac{\partial \u}{\partial t} + v|\Grad\u|=0 \\ [2mm]
\displaystyle v=-\frac{K[\u]}{\rho_s} \big(C_s-C_b\big)
\\ [4mm]
\displaystyle \frac{d C_b}{dt}=\frac{C_s-C_b}{\V}\int_{\{\phi=0\}}K ds
\end{array}
\right. ,
\quad C_b<C_s,
\end{equation}
with initial conditions given by \eqref{phi0} and $C_b(0)=0$.
When $C_b>C_s$ dissolution is stopped and recrystallization in the bulk fluid occurs:
$$
\left\{
\begin{array}{l}
\displaystyle\frac{\partial \u}{\partial t}=0    \\ [4mm]
\displaystyle \frac{d C_b}{dt}=\big(C_s-C_b\big)k_{rb}, 
\end{array}
\right. ,
\quad C_b\geq C_s.
$$
If the particle comes to a complete dissolution, we can directly set $K=0$ in \eqref{mainequation_singleparticle}.

\subsubsection*{Multiple particles}
In the case of $N_P>1$ particles dissolving in the same volume, one should track separately every particle. 
Differently from classical fields of application of the level-set method, here particles do not interact to each other: their dynamics are independent, but they all contribute to the concentration $C_b$ by summing all single contributions. 
Therefore, in the presence of multiple particles the equation \eqref{mainequation_singleparticle} is simply generalized by
\begin{equation}\label{mainequation_multipleparticle}
\left\{
\begin{array}{l}
\displaystyle\frac{\partial \u_p}{\partial t}+v_p|\Grad\u_p|=0 \\ [2mm]
\displaystyle v_p=-\displaystyle\frac{K_p[\u_p]}{\rho_s} \big(C_s-C_b\big)  \\ [4mm]
\displaystyle\frac{d C_b}{dt}=\frac{C_s-C_b}{\V}\sum_{p=1}^{N_P}\int_{\{\u_p=0\}}K_p ds
\end{array}
\right.,\qquad p=1,\ldots, N_P.
\end{equation}
\bigskip

If we have $M>1$ \emph{identical} particles, one can avoid to track them separately. Indeed, their contribution to $C_b$ will be obtained by 
simply multiplying the term $\int_{\{\u_p=0\}}K_p ds$ by $M$.

\section{Numerical simulations}\label{sec:simulations}
In this section we report the results of some numerical tests obtained by using the level-set methodology.
As in \cite{abrami2022EJPB}, we prefer to work with a dimension-reduced 2D dynamics. 
This greatly simplifies computations still allowing us to investigate qualitatively the impact of the particle shape and model's parameters on the dissolution dynamics.
Clearly, all considerations and results related to the particle's perimeter (resp.\ area) in 2D will be naturally transferred to its surface (resp.\ volume) in 3D.

In Sect.\ \ref{sec:numerics_SP}, we consider the case of a single particle, for three drugs: theophylline (25$^{\circ}$C and 37$^\circ$C), griseofulvin (37$^\circ$C), and nimesulide (37$^\circ$C), with parameters as in Table \ref{tab:parameters}. 
In particular, theophylline (25$^{\circ}$C) is considered to facilitate the comparison with results reported in \cite{abrami2022EJPB}. 
These drugs are widely used in the pharmaceutical field \cite{abrami2020ADMET} and they are representative of 
soluble and wettable drugs (theophylline), 
of poorly soluble but wettable drugs (griseofulvin), 
and of poorly soluble and poorly wettable drugs (nimesulide).

In Sect.\ \ref{sec:numerics_MP}, instead, the more realistic polydisperse case is analyzed, for griseofulvin (37$^\circ$C) only. 

All numerical details are postponed in Appendix \ref{app:numerics}.

\begin{table}[h!]
    \centering
    \rotatebox{90}{
    \begin{minipage}{1.0\textheight}
      \centering
    \begin{tabular}{l|l|l|l|l|l|l}
         Parameter & Symbol & Unit & Theoph.\ 25$^\circ$C & Theoph.\ 37$^\circ$C & Griseof.\ 37$^\circ$C & Nimes.\ 37$^\circ$C \\ \hline
         Solid drug density                  & $\rho_s$      & kg/m$^3$ & 1490                 & 1490                 & 1495                  & 1476 \\
         Initial drug solubility             & $C_{s0}$      & kg/m$^3$ & 11.6                 & 12.495               & 0.494                 & 4.108\\
         Final drug solubility               & $C_{sf}$      & kg/m$^3$ & 6.1                  & 6.569                & 0.025                 & 0.028\\
         Surface recrystallization const.    & $k_r$         & s$^{-1}$ & $6\cdot 10^{-3}$     & $6\cdot 10^{-3}$     & $8.8\cdot 10^{-3}$    & $1.3\cdot 10^{-2}$\\
         Bulk recrystallization const.       & $k_{rb}$      & s$^{-1}$ & $6.6\cdot 10^{-3}$   & $5.7\cdot 10^{-3}$   & $8.36\cdot 10^{-3}$   & $1.235\cdot 10^{-2}$\\
         Interfacial flat mass transfer coeff. & $k_m^\infty$  & m/s      & $3.7\cdot 10^{-3}$   & $2\cdot 10^{-3}$     & $0.126$               & $1.8\cdot 10^{-7}$  \\
         Drug diffusivity in water           & $D$           & m$^2$/s  & $6.2\cdot 10^{-10}$  & $8.2\cdot 10^{-10}$  & $7.057\cdot 10^{-10}$ & $7.388 \cdot 10^{-10}$\\
         Fluid density                       & $\rho_f$      & kg/m$^3$ & 1000                 & 993                  & 993                   & 993\\
         Fluid dynamic viscosity             & $\eta_f$      & Pa s     & $10^{-3}$            & $6.91 \cdot 10^{-4}$ & $6.91 \cdot 10^{-4}$  & $6.91 \cdot 10^{-4}$\\
         Fluid kinematic viscosity           & $\nu_f$       & m$^2$/s  & $10^{-6}$            & $6.96 \cdot 10^{-7}$ & $6.96 \cdot 10^{-7}$  & $6.96 \cdot 10^{-7}$\\
         Parameter in $k_m$                  & $\alpha$      & m$^3$    & $10^{-15}$           & $10^{-15}$           & $10^{-15}$            & $10^{-15}$\\
         Tolman length                      & $d_T $        & m        & $10^{-9}$            & $2.6 \cdot 10^{-10}$ & $3.10 \cdot 10^{-10}$ & $2.82 \cdot 10^{-10}$
    \end{tabular}
    \caption{Physical parameters of the three considered drugs used in numerical simulations.}
    \label{tab:parameters}
    \end{minipage} 
    }
\end{table}

\clearpage
\subsection{Dissolution of a single particle}\label{sec:numerics_SP}
\subsubsection*{Test 1a: Circular particle (theophylline 25 $^\circ$C)}
In this preliminary test we consider a single circular particle of theophylline 25$^\circ$C with initial radius $R_0=50$ $\mu$m. 
The main goal here is to compare the result obtained by the level-set method with that reported in \cite{abrami2022EJPB}, as well as with the solution of the simplified problem detailed in Rem.\ \ref{rem:cerchi}. 
Eq.\ \eqref{mainequation_singlecircularparticle} is solved using the third-order Runge-Kutta method and can be assumed to be the `exact' solution.

\begin{figure}[h!]
    \centering
    \includegraphics[width=0.49\textwidth]{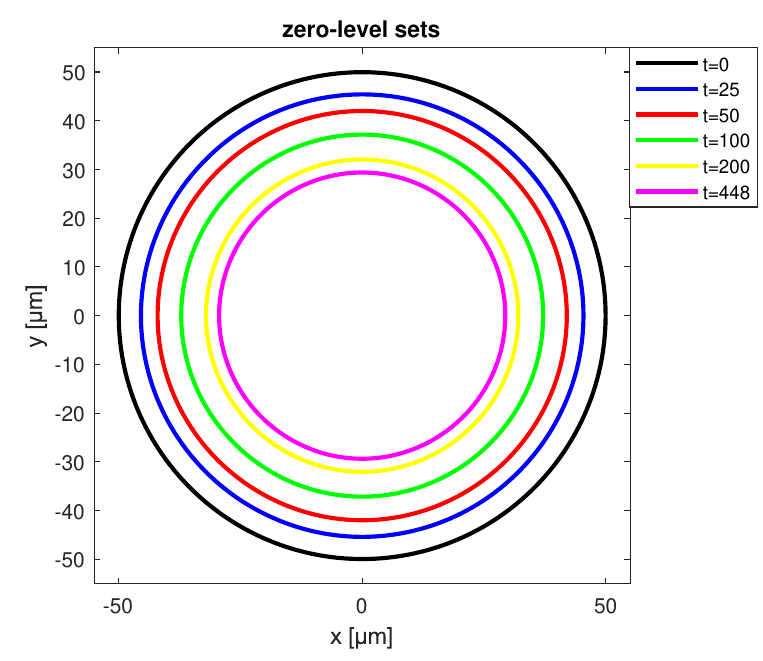}
    \includegraphics[width=0.49\textwidth]{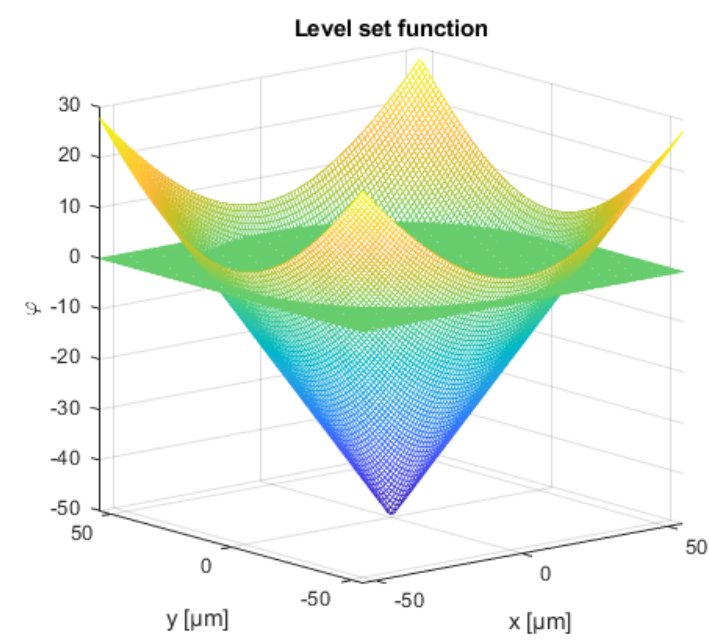}
    \caption{Test 1a. Evolution of the particle with $V^+=150$ (left) and its level-set function $\u$ (right).}
    \label{fig:T1a-evoluzione}
\end{figure}

Fig.\ \ref{fig:T1a-evoluzione} shows the evolution of the particle with $V^+=150$ and its level-set function $\u$, cut at the zero level set to get the particle's profile at a given time. 
In this case the recrystallization process is present. 
No appreciable difference can be seen comparing the level-set solution with the results obtained in \cite{abrami2022EJPB}. 

\begin{figure}[h!]
    \centering
    \includegraphics[width=0.55\textwidth]{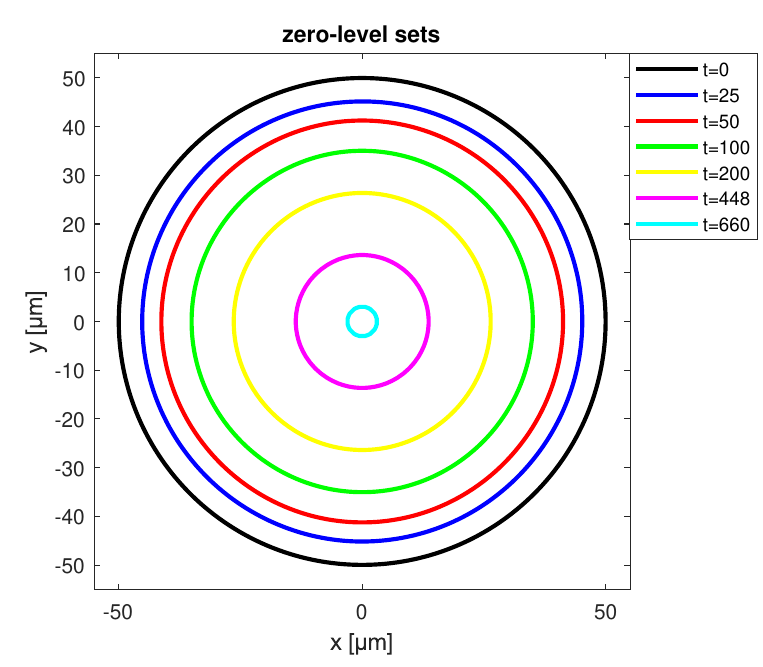}
    \includegraphics[width=0.44\textwidth]{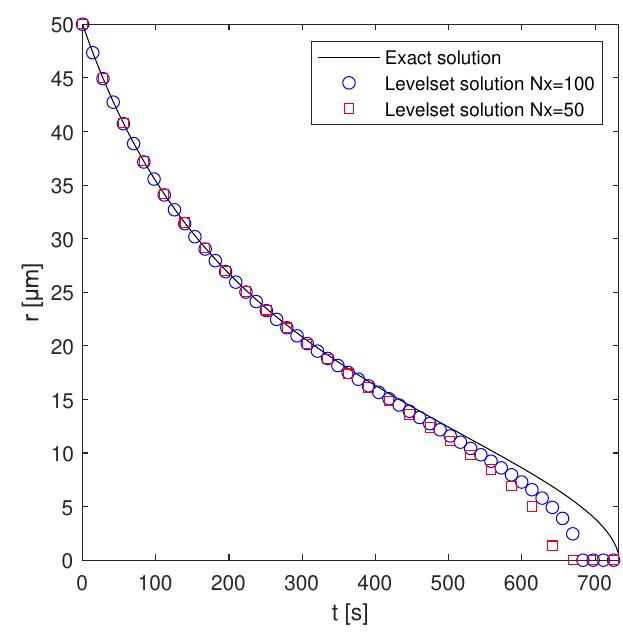}
    \caption{Test 1a. Evolution of the particle with $V^+=300$ (left) and the evolution of the particle's radius $R(t)$ on a coarse/fine numerical grid (right).}
    \label{fig:T1a-phi_e_errore}
\end{figure}

Fig.\ \ref{fig:T1a-phi_e_errore} reports the evolution of the particle with $V^+=300$ and the evolution of the particle's radius $R(t)$, comparing the exact solution and the level-set solution for two space steps $\Delta x$'s (keeping the CFL ratio constant, see Appendix A).
In this case the particle reaches the complete dissolution. Again, no significant difference is seen comparing the level-set solution with the results obtained in \cite{abrami2022EJPB}. 
On the other hand, the level-set method does not catch very well the evolution of the particle when it is very small. 
This is due to the fact that the particle size becomes comparable with that of the grid cell. 
As expected, increasing the number of grid points (i.e.\ reducing $\Delta x$) the error diminishes. 

\begin{figure}[h!]
    \centering
    \includegraphics[width=0.6\textwidth]{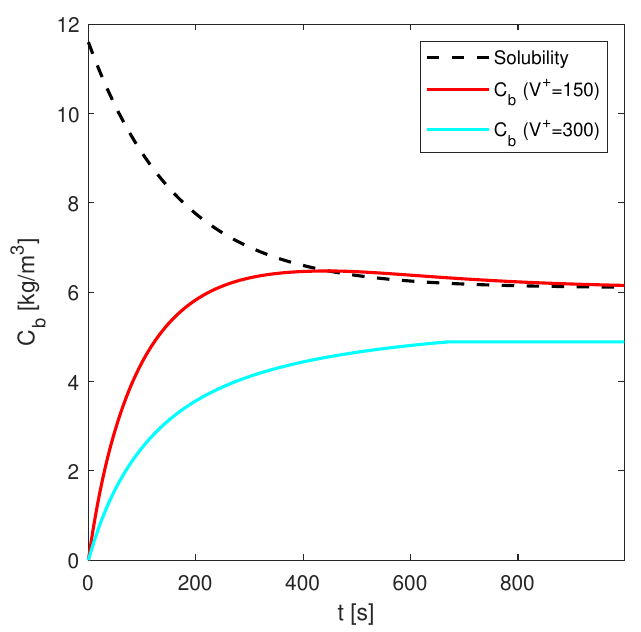}
    \caption{Test 1a. Evolution of $C_b$ for $V^+=150$ and $V^+=300$.}
    \label{fig:T1a-Cb}
\end{figure}

Fig.\ \ref{fig:T1a-Cb} shows the curves of $C_b$ for $V^+=150$ and $V^+=300$.
It is interesting to underline that when the volume of the dissolution medium is large ($V^+=300$) recrystallization does not occur as drug concentration (blue line) never exceeds the time dependent drug solubility (dotted line). On the contrary, this event takes place for smaller $V^+$.

\subsubsection*{Test 1b: Circular particle (theophylline 37$^\circ$C)}
To consider the effect of temperature, we consider now a single circular particle of theophylline 37$^\circ$C with initial radius $R_0=50$ $\mu$m. 

\begin{figure}[h!]
    \centering
    \includegraphics[width=0.49\textwidth]{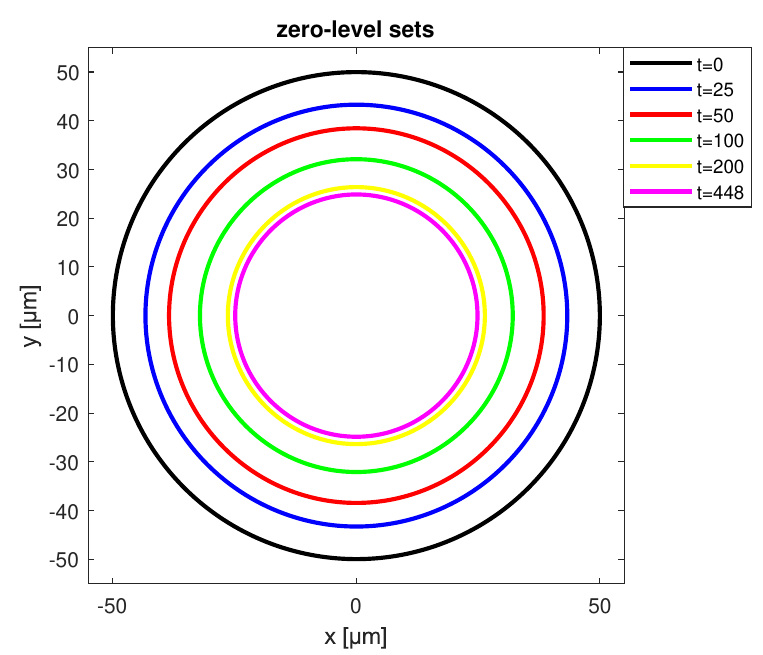}
    \includegraphics[width=0.49\textwidth]{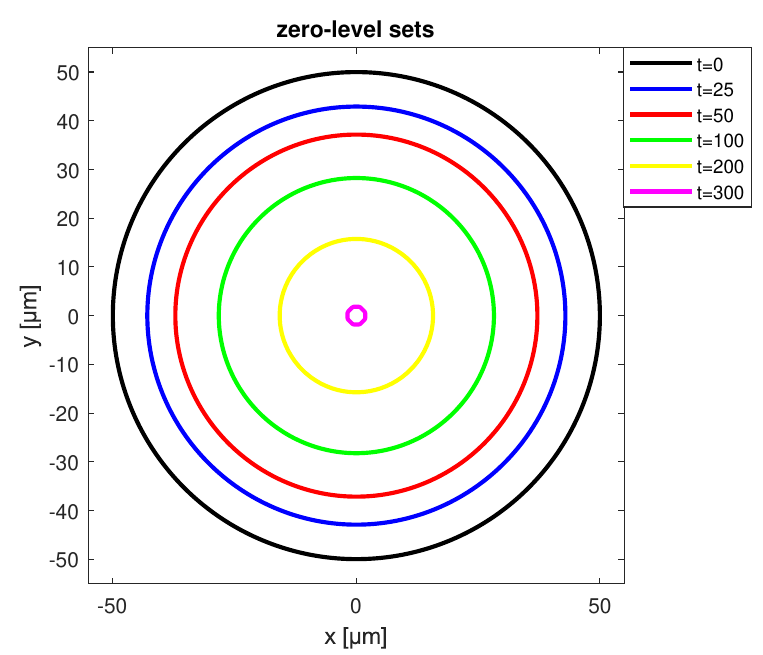}
    \caption{Test 1b. Evolution of the particle with $V^+=150$ (left) and $V^+=300$ (right).}
    \label{fig:T1b-evoluzione}
\end{figure}

\begin{figure}[h!]
    \centering
    \includegraphics[width=0.6\textwidth]{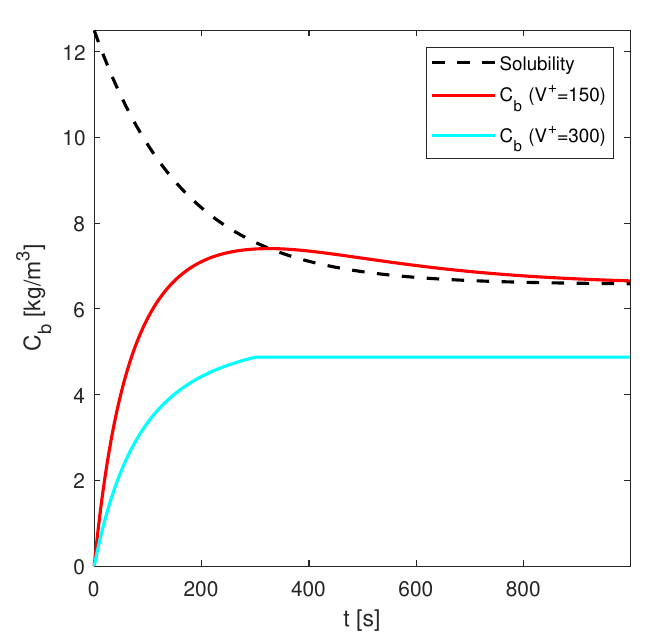}
    \caption{Test 1b. Evolution of $C_b$ for $V^+=150$ and $V^+=300$.}
    \label{fig:T1b-Cb}
\end{figure}

Fig.\ \ref{fig:T1b-evoluzione} (resp.\ \ref{fig:T1b-Cb}) shows the evolution of the particle surface (resp. $C_b$) with $V^+=150$ and $V^+=300$.
Comparison between Tests 1a and 1b reveals how the increase of temperature enhances the dissolution kinetics. 
This is due to the increased theophylline diffusion coefficient and solubility in the dissolution medium (water).

\subsubsection*{Test 2: From circle to square (theophylline 37$^\circ$C)}
In this test we investigate the impact of the shape of a particle of theophylline on its evolution. We use superellipses (actually supercircles in this case), with exponent $n=2,3,39$ (see Appendix \ref{app:inizialization} for details) in order to smoothly move from a circle to a square \emph{with same initial areas} ($R_0=50$ $\mu$m for the circle and $V^+=150$).

\begin{figure}[h!]
    \centering
    \includegraphics[width=0.49\textwidth]{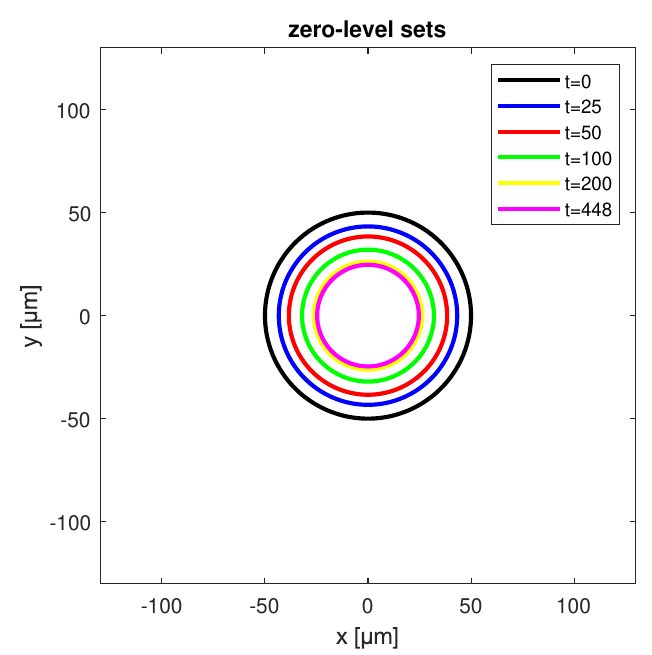}
    \includegraphics[width=0.49\textwidth]{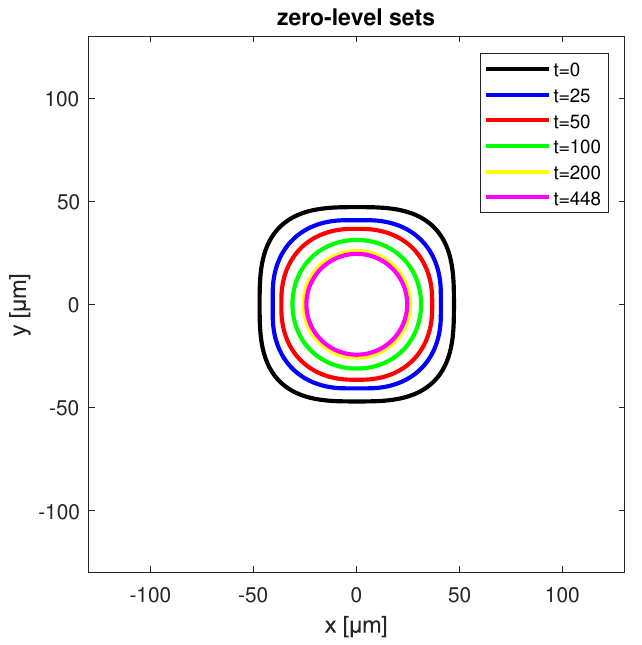} \\
    \includegraphics[width=0.49\textwidth]{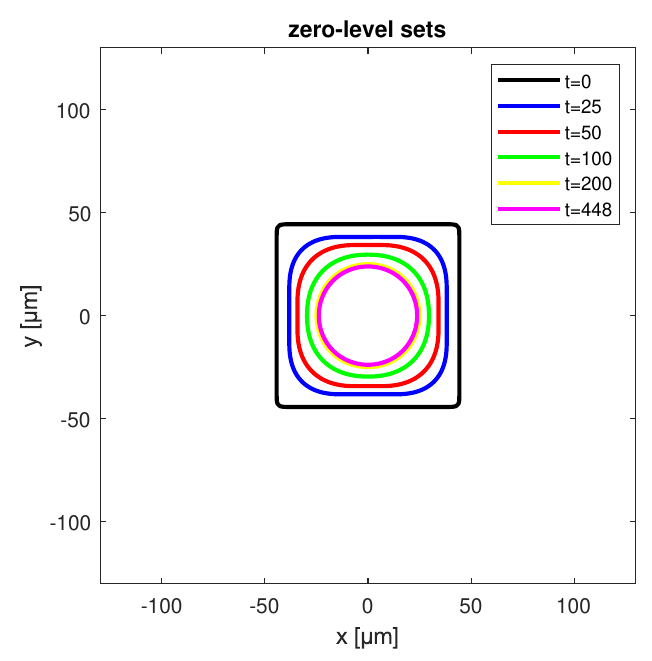}
    \includegraphics[width=0.49\textwidth]{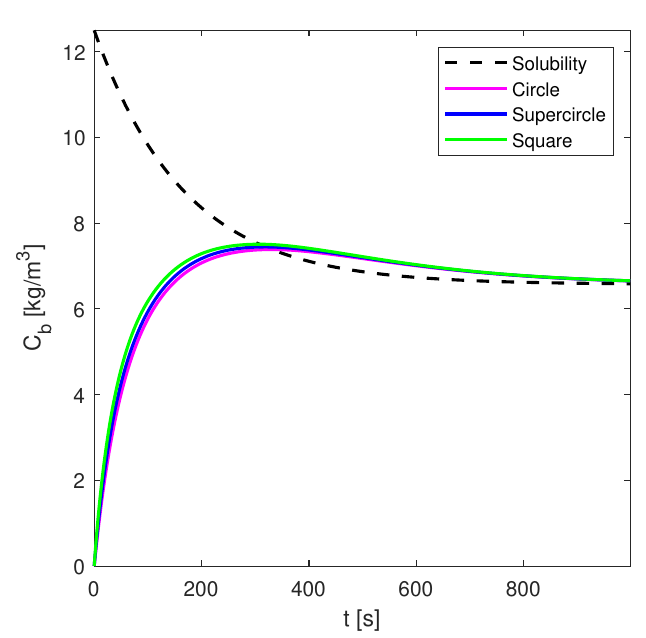}
    \caption{Test 2. Dissolution profiles of particles with different initial shapes. 
    Top-left: circle, 
    top-right: supercircle ($n=3$), 
    bottom-left: square (actually a supercircle with $n=39$). 
    Bottom-right: evolution of $C_b$. 
    }
    \label{fig:T2}
\end{figure}

\begin{figure}[h!]
    \centering
    \includegraphics[width=0.6\textwidth]{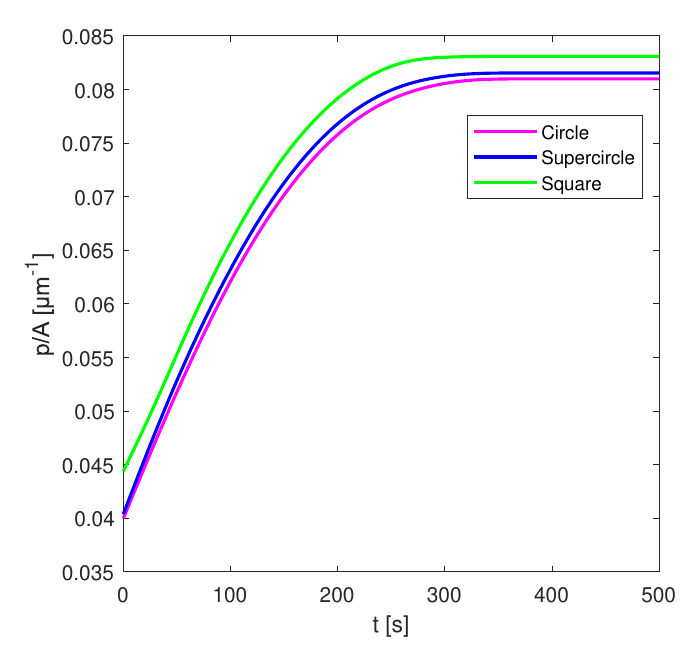}
    \caption{Test 2. Evolution of $p/A$ for the three shapes.}
    \label{fig:T2-ratio-p/A}
\end{figure}

Fig.\ \ref{fig:T2} shows the dissolution profiles and the comparison of $C_b$'s. We note that the  sharp edges turn to be become rounded, and the square tend to become a circle. Consequently the evolution of $C_b$ is similar in the three cases.
Both perimeter $p$ and area $A$ diminish in time. Their ratio $p/A$ is also not constant, although the same trend for the three shapes is maintained, see Fig.\ \ref{fig:T2-ratio-p/A}. 
The highest values are reached by the square and the lowest by the circle.

\begin{figure}[h!]
    \centering
    \includegraphics[width=0.49\textwidth]{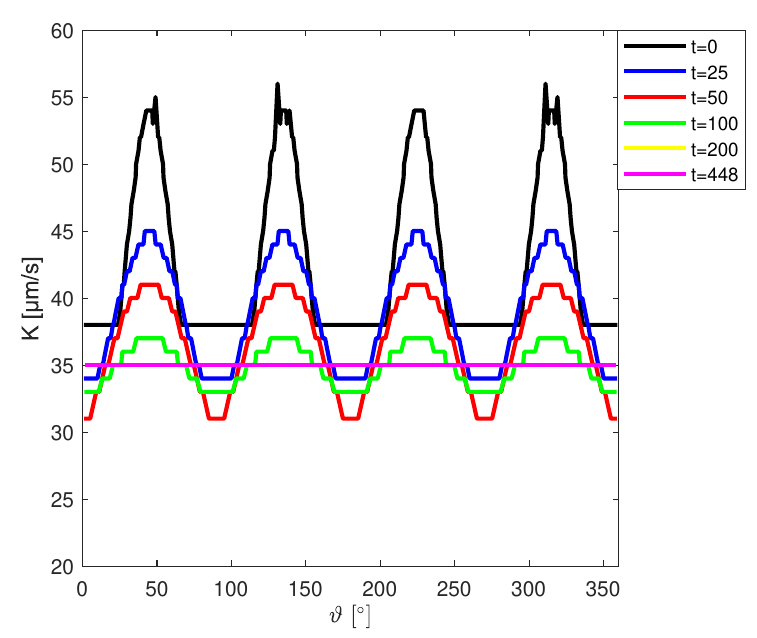}
    \includegraphics[width=0.49\textwidth]{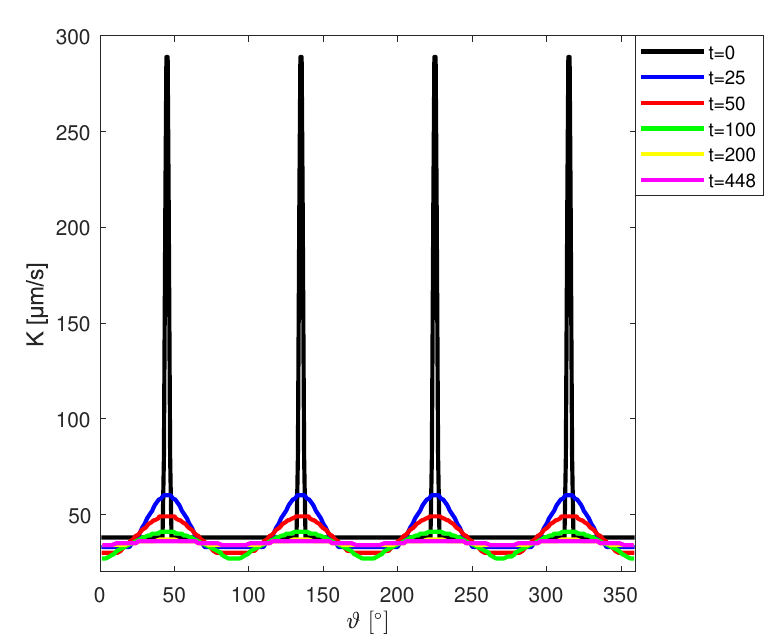}
    \caption{Test 2. Mass transfer coefficient $K$ as a function of the angle $\theta\in[0,360^\circ]$, at various times, for the supercircle ($n=3$) (left) and the square (right).}
    \label{fig:T2_K}
\end{figure}

Finally, Fig.\ \ref{fig:T2_K} shows the mass transfer coefficient $K$ at any point of the particle boundary, for the circle and the square.
One can observe that $K$ largely increases at corners, cf.\ \cite{abrami2022EJPB}. 
Note that stair-like behavior of the function is just a grid effect due to the numerical discretization.

\subsubsection*{Test 3: From square to rectangle (theophylline 37$^\circ$C)}
In this test we investigate further the impact of the shape of a particle of theophylline on its evolution using again superellipses with  exponent $n=39$, smoothly moving from a square to a rectangle of increasing aspect ratio, \emph{with same initial areas} and $V^+=150$.

\begin{figure}[h!]
    \centering
    \includegraphics[width=0.49\textwidth]{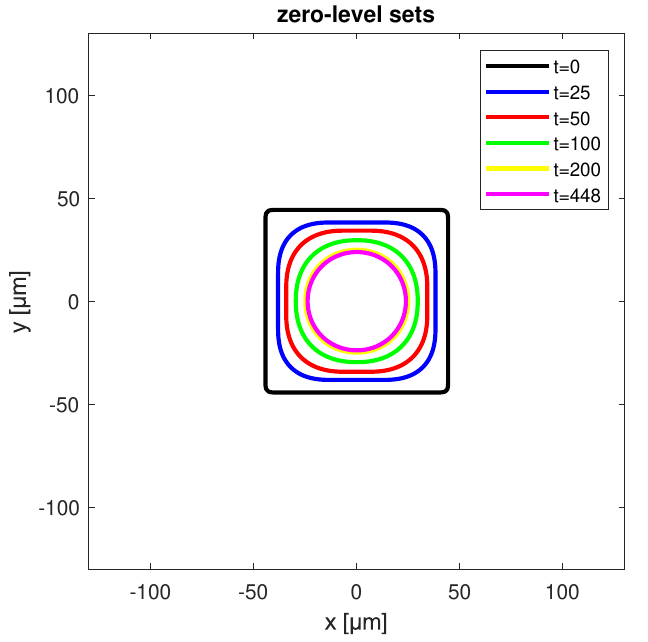}
    \includegraphics[width=0.49\textwidth]{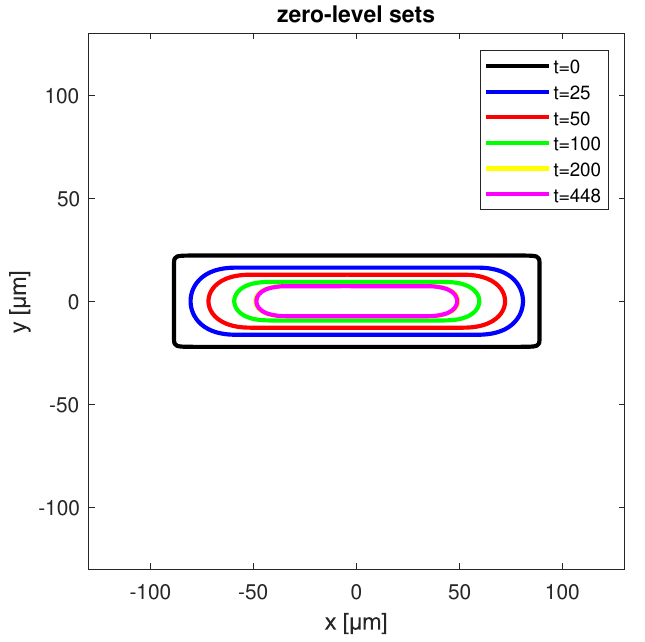} \\
    \includegraphics[width=0.49\textwidth]{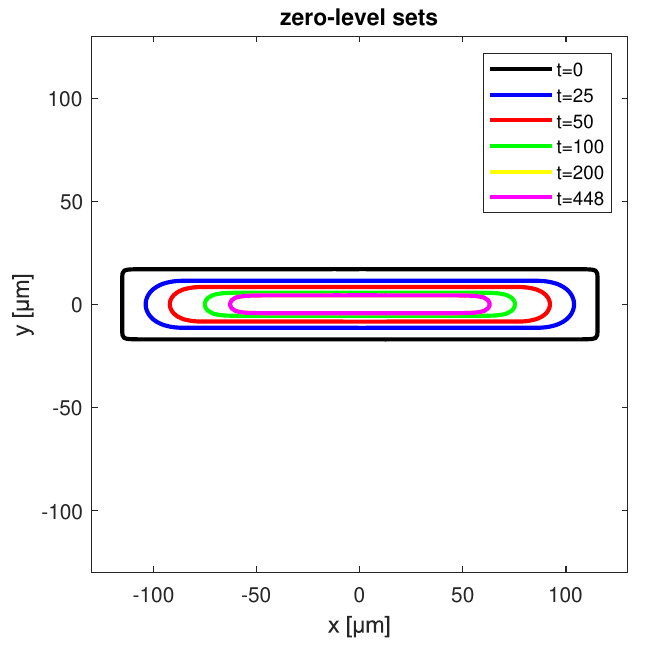}
    \includegraphics[width=0.49\textwidth]{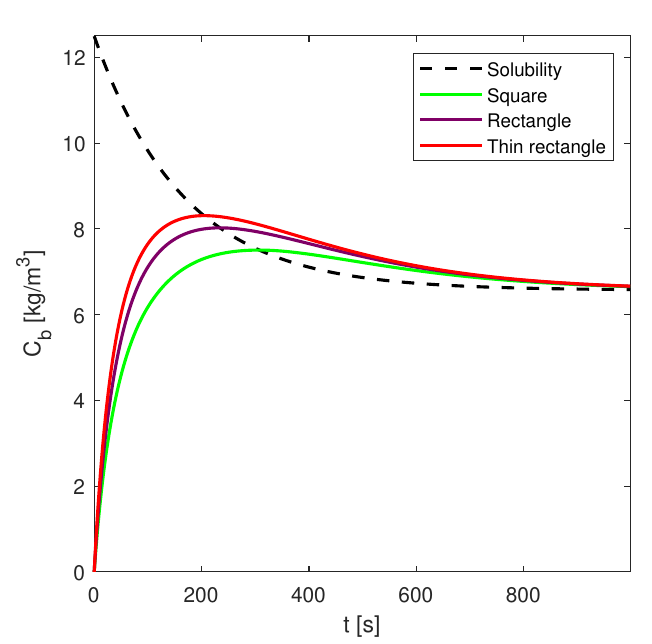}
    \caption{Test 3. Dissolution profiles and comparison of $C_b$, for 3 rectangles with increasing aspect ratio. Initial dimensions are
    87 $\mu$m $\times$ 87 $\mu$m (top-left), 
    177 $\mu$m $\times$ 44 $\mu$m (top-right), 
    231 $\mu$m $\times$ 34 $\mu$m (bottom-left).}
    \label{fig:T3}
\end{figure}

Fig.\ \ref{fig:T3} shows the dissolution profiles and the comparison of $C_b$'s. 
This time, the evolutions of $C_b$'s are not similar, meaning that impact of the shape is not negligible. 
%

\begin{figure}[h!]
    \centering
    \includegraphics[width=0.6\textwidth]{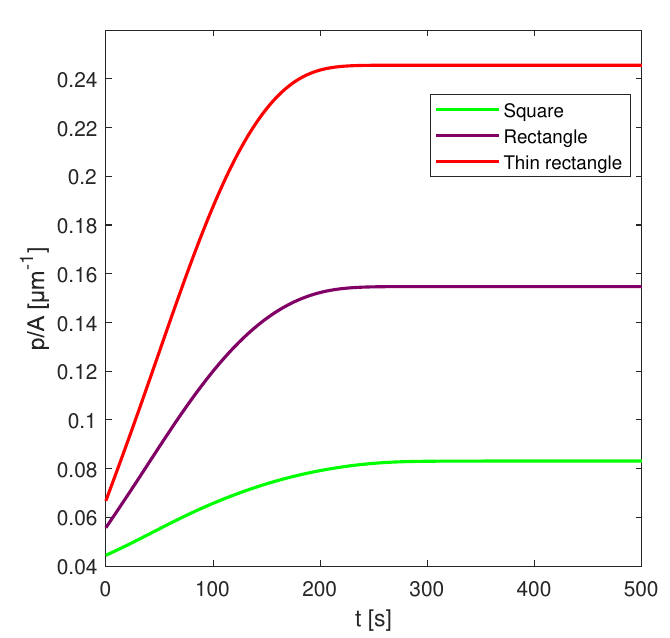}
    \caption{Test 3. Evolution of the ratio $p/A$ for the three shapes.}
    \label{fig:T3-ratio-p/A}
\end{figure}

Fig.\ \ref{fig:T3-ratio-p/A} shows the evolution of the ratio $p/A$ for the three shapes. 
The highest value is reached by the elongated rectangle and the lowest by the square. 
Here the difference among shapes is more significant than that observed in Fig. \ref{fig:T2-ratio-p/A}. 
This suggests that the evolution of $C_b$ is influenced especially by the ratio $p/A$, rather than by the presence of corners.
This evidence is very important from a practical point of view, as the reproducibility of dissolution experiments mainly depends on the particles surface/volume ratio rather than on the presence of possible irregularities (corners) on particles surface.

\subsubsection*{Test 4: Comparison between three drugs}
Here we consider three circular particles of theophylline, griseofulvin and nimesulide (37$^\circ$C), all of them with initial radius $R_0=250$ $\mu$m and $V^+=150$. Even if such a radius is not physically feasible for all the three drugs, we perform this test to compare them all.
Moreover, the value of $C_b$ is normalized over $C_{sf}$ for an effective comparison.
Beside the comparison of the dissolution of the three drugs, we also aim at investigating the role of the parameter $\alpha$ appearing in \eqref{km} and ruling the interfacial mass transfer coefficient $k_m$, which is very difficult to measure experimentally.
Therefore, we use \textit{in silico} experiments to understand if it is possible to have a clue on $\alpha$ observing the particle's evolution \textit{a posteriori}. 
To this goal, we use two values for $\alpha$, i.e.\ $\alpha_1=10^{-15}$ and $\alpha_2=10^{-2}$.

\begin{figure}[h!]
    \centering
    \includegraphics[width=0.49\textwidth]{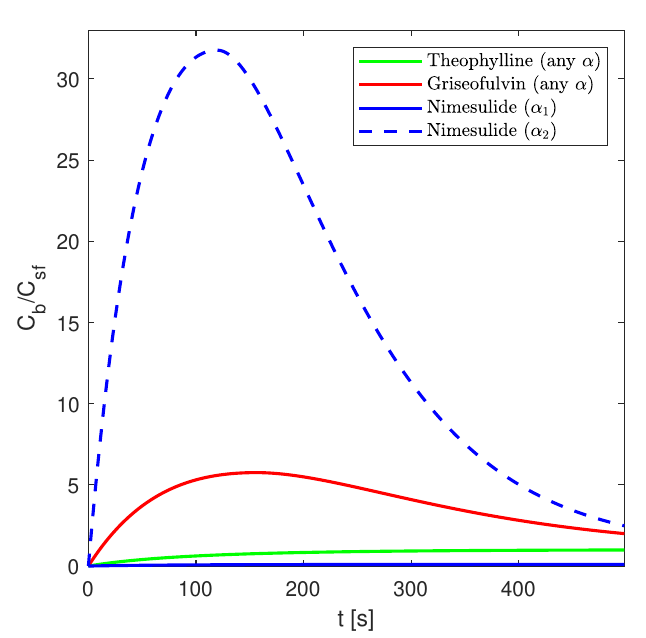}
    \includegraphics[width=0.49\textwidth]{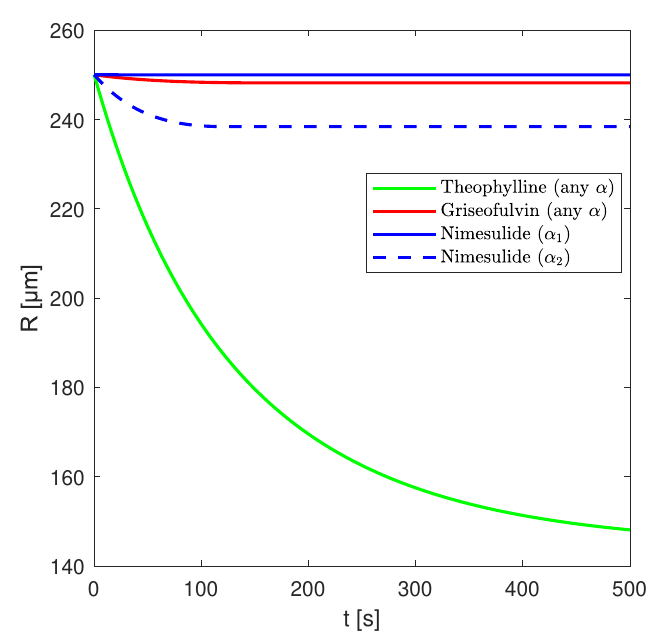}
    \caption{Test 4. Evolution of $C_b/C_{sf}$  (left) and evolution of the radii $R$ (right) of circular particles of theophylline, griseofulvin, and nimesulide, with $\alpha=\alpha_1,\alpha_2$.}
    \label{fig:T4}
\end{figure}

Fig.\ \ref{fig:T4} shows, for the three drugs, the evolution of $C_b/C_{sf}$ and the evolution of the radii, with $\alpha=\alpha_1, \alpha_2$. 
It is observed that, among the three drugs, nimesulide is the only one for which $\alpha$ has a relevant impact on dissolution.
This is not surprising as nimesulide is the only drug, among those considered, showing poor wettability, that is connected to an important role played by $k_m$ on the dissolution kinetics.

\clearpage 
\subsection{Dissolution of multiple particles}\label{sec:numerics_MP}
In the pharmaceutical experimental practise, drug powders and compounds appear smashed in 
a mixture of grains of different sizes and shapes, before being dissolved.
In this section we analyze the effect of a polydisperse mixture of micro-sized particles in a variety of sizes and shapes.

In the following numerical tests we have chosen to work with griseofulvin (37$^\circ$ C) because its particles' size range of a typical mixture is smaller (1-15 $\mu$m) \cite{mosharraf1995, naseem2004} than those of the other considered drugs. 
This simplifies the numerical treatment because we can use the same numerical grid for all particles' sizes.

In the following we still denote by $p$ and $A$ the perimeter and the area of the particles, respectively, here meaning the \emph{total} perimeter and the \emph{total} area.

\subsubsection*{Test 5: Effect of the particles' size}
\begin{figure}[h!]
    \centering
    \includegraphics[width=0.6\textwidth]{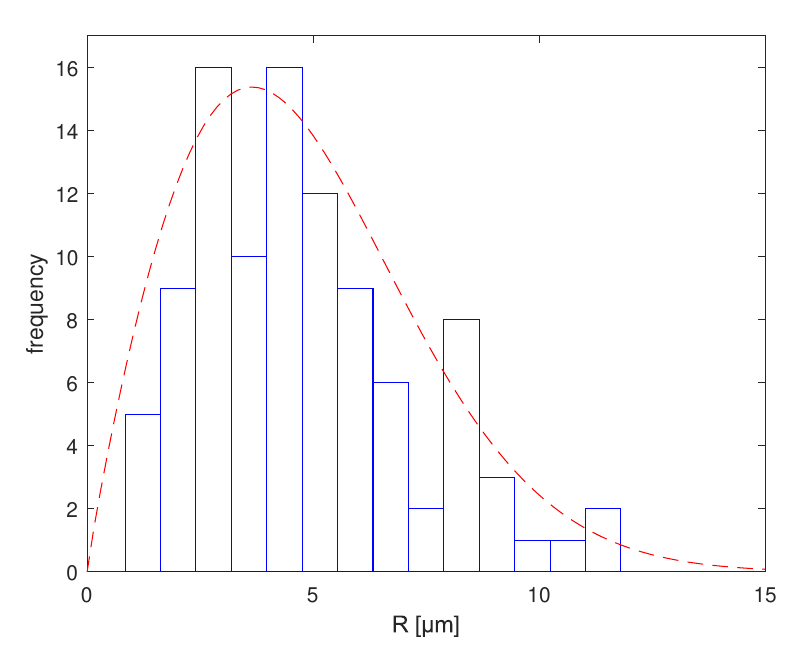}
    \caption{Test 5. Particles' size (radius) distribution for 100 different circles (case c).
    Weibull distribution parameters: $\lambda=5.4$, $k=1.9$, $x_0=0$ \cite[Fig.\ 1(a)]{abrami2020ADMET}), see also Appendix \ref{app:weibull}.
    }
    \label{fig:M-T5-distribuzioneraggi}
\end{figure}

In this test we compare the dissolution of three families of polydisperse drug particles sharing the same area ($A=9,125.78$ $\mu$m$^2$), but having different perimeter $p$:  
\begin{itemize}
    \item[a)] a single circular particle with initial radius $R_0=53.90$ $\mu$m ($p=338.64$ $\mu $m);
\item[b)] 100 equal circular particles of initial radius $R_0=5.39$ $\mu$m ($p=3,386.41$ $\mu $m);
\item[c)] 100 circular particles with different initial radii, randomly chosen according to the Weibull distribution, see Fig.\ \ref{fig:M-T5-distribuzioneraggi} and Appendix \ref{app:weibull} ($p=3,019.70$ $\mu $m). 
\end{itemize}
In all cases $V^+=1,000$.  
\begin{figure}[h!]
    \centering
    \includegraphics[width=0.49\textwidth]{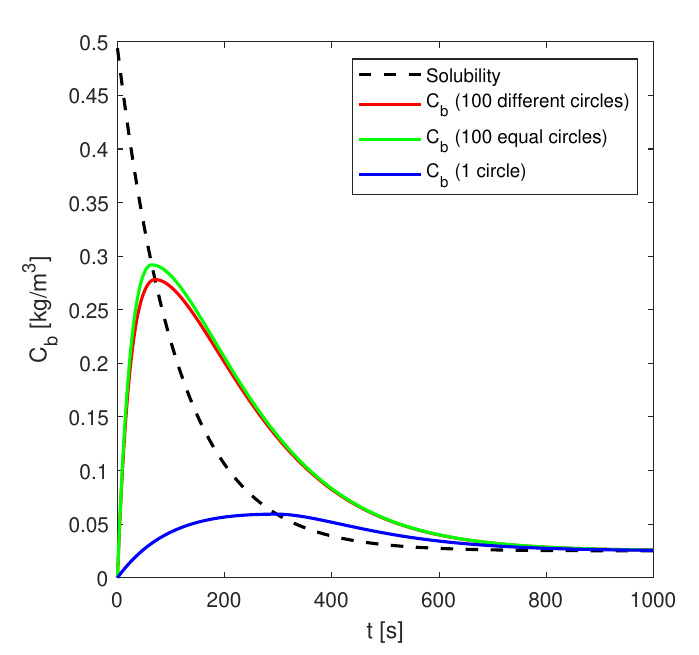}
    \includegraphics[width=0.47\textwidth]{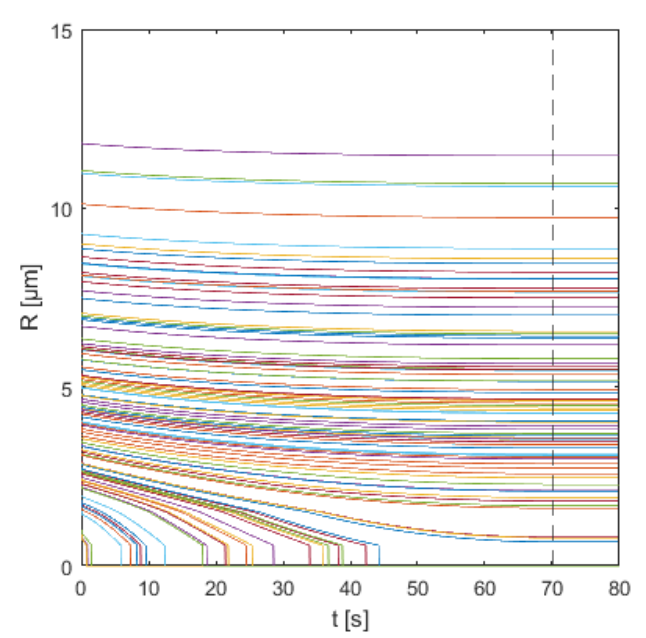}
    \caption{Test 5. Evolution of $C_b$ for cases a, b, c (left), and evolution of  particles' radii for case c only (right). 
    In case c, the 24 smallest particles reach full dissolution. 
    }
    \label{fig:M-T1}
\end{figure}
Fig.\ \ref{fig:M-T1} shows that the single large particle (case a) and the 100 equal particles (case b) do not reach full dissolution, while the 24 smallest particles of case c undergo a full dissolution. 
As expected, $C_b$ grows more rapidly in the case of many small particles, but the difference between cases b and c is mainly due to the different number of fully dissolved particles, rather than to the difference of $p/A$, cf.\ next test.

\clearpage 
\subsubsection*{Test 6: Effect of the particles' shape}
\begin{figure}[h!]
    \centering
    \includegraphics[width=0.5\textwidth]{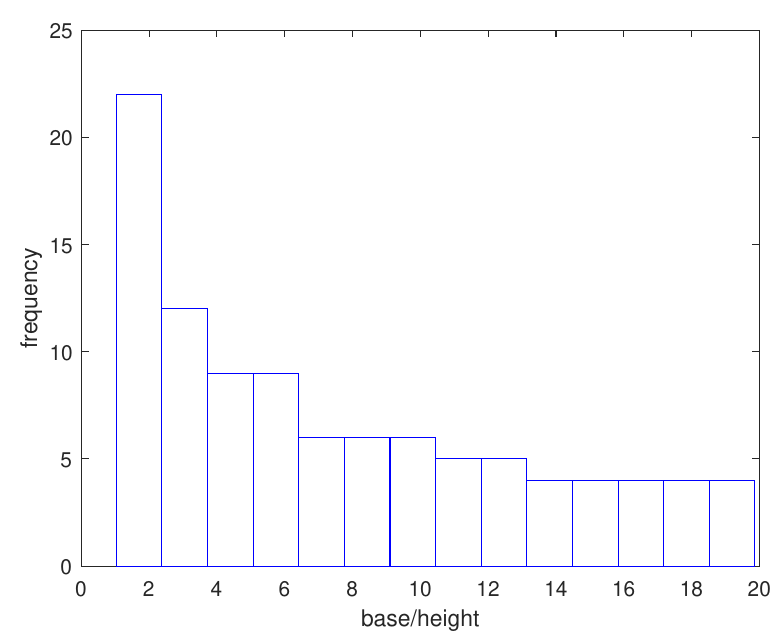}
    \caption{Test 6. Rectangular particles' aspect ratio distribution.}
    \label{fig:M-T2-distribuzioneratios}
\end{figure}

\begin{figure}[h!]
    \centering
    \includegraphics[width=0.49\textwidth]{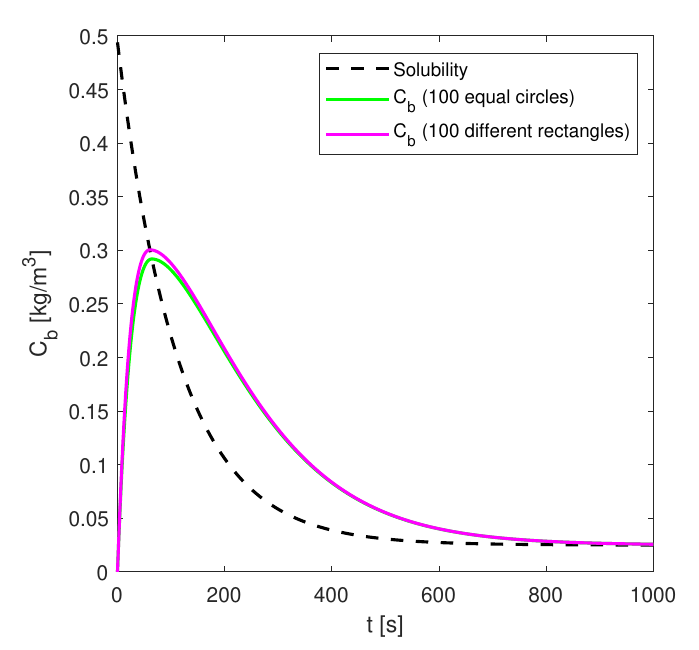}
    \includegraphics[width=0.46\textwidth]{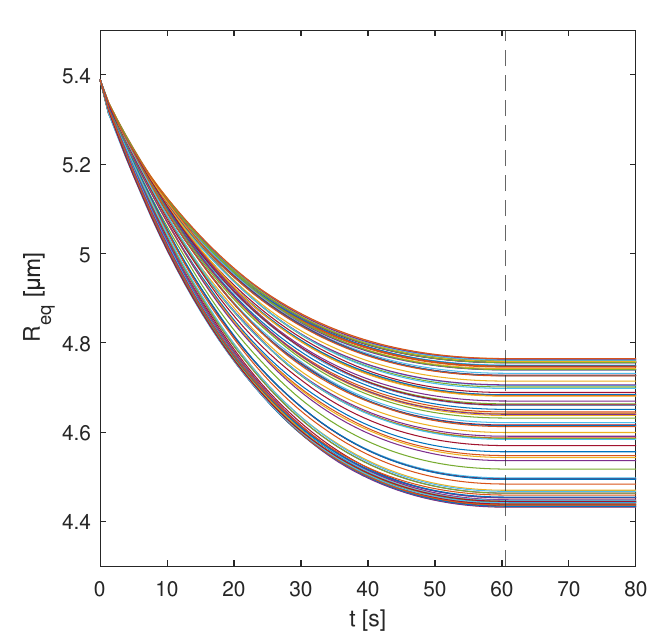}
    \caption{Test 6. Left: Evolution of $C_b$ for 100 equal circles (case b) and 100 rectangles (case d). Right: evolution of the 100 rectangles' equivalent radii. 
    No particle reaches full dissolution.
    }
    \label{fig:M-T2}
\end{figure}

To emphasize the influence of the shape, here  two families of 100 particles, with the only difference in shape, are considered. 
The area of particles is the same, both individually as well as the total one.
In particular, we compare the dissolution of  100 equal circular particles (case b of the previous test) with 
\begin{itemize}
    \item[d)] 100 rectangular particles with different aspect ratios, with $p=5,805.18$ $\mu$m, distributed as in Fig.\ \ref{fig:M-T2-distribuzioneratios}. 
\end{itemize}
We choose again $V^+=1,000$. 
Fig.\ \ref{fig:M-T2} shows the evolution of $C_b$ for case b and d, and the evolution of the rectangles' equivalent radii. 
Note that in this case no particle attains full dissolution.
One can see here that the shape has a minor impact on the dissolution. 
Actually this is not in contrast with results of Test 3, for two reasons: griseofulvin is less sensitive to $p/A$ and $V^+$ is relatively small to bring out a significant difference: for larger values of $V^+$ a more remarkable effect is indeed visible.

\subsubsection*{Test 7: The polydisperse case}
To mimic a more realistic situation of a mixture with a variety of shape and sizes, in this final test the simultaneous dissolution of 20 circles, 50 ellipses, and 30 rectangles with randomly chosen size and aspect ratios, is simulated. We have $p=6,072.92$ $\mu$m and $A=48,456.72$ $\mu$m$^2$. 
The result is also compared with that of 61 equal circles of radius 15.96 $\mu$m, which share the same $p$ and the same $A$ with the mixture.
Also, $V^+=10,\!000$.

\begin{figure}[h!]
    \centering
    \includegraphics[width=0.49\textwidth]{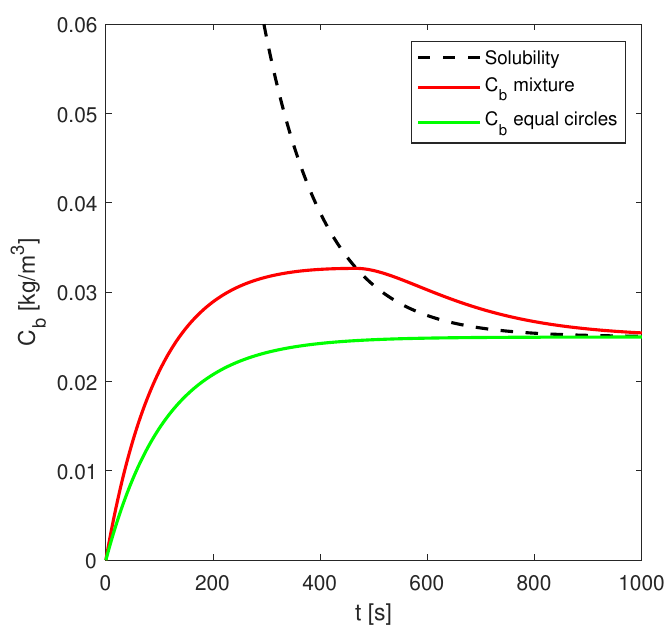}
    \includegraphics[width=0.49\textwidth]{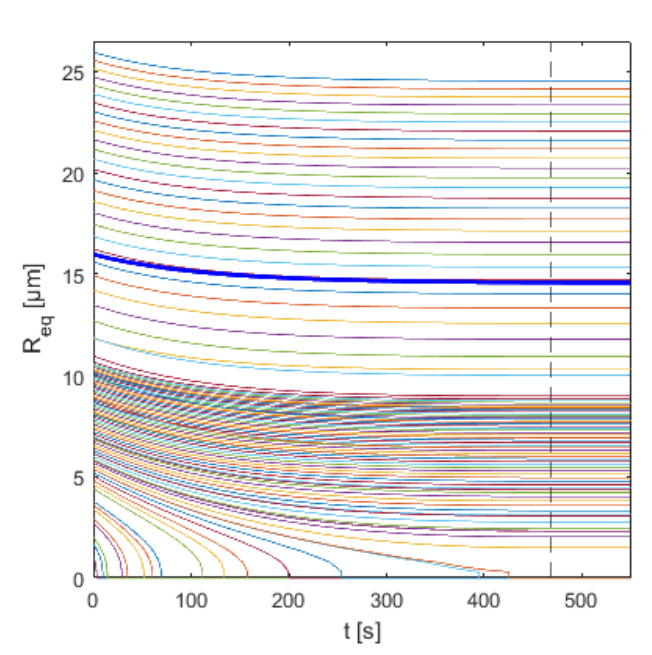}
    \caption{Test 7. Left: evolution of $C_b$ for 100 particles with different size and shape and 100 circles, same $p$ and $A$. 
    Right: evolution of the equivalent radii of the 100 particles with various shape an size (colored thin curves, 15 particles reach a full dissolution), and evolution of the radius of one of the 61 circular particles (blue thick curve).}
    \label{fig:M-T4}
\end{figure}

Fig.\ \ref{fig:M-T4} shows the evolution of $C_b$ and that of the particles' equivalent radii.
Results make it clear that particles shape and distribution (inside each shape) reflects in a small but not negligible variation of the dissolution profile.

\section{Conclusions}
In order to understand the role played by the different concurrent phenomena taking place upon the dissolution of an ensemble of different particles (also referred to as DRT (Dissolution Rate Test \cite{abrami2020ADMET}), this paper focuses on the realization of a mathematical model able to consider the majority of them, with particular care devoted to particles shape, size distribution and surface properties. 
To this purpose, the three drugs (theophylline, griseofulvin and nimesulide), differing for shape, size distribution and for what concerns their water solubility and wettability, have been considered for the simulations of DRT. 
The numerical solution of the proposed model has been achieved by the level-set method \cite{osher1988JCP}, a reliable approach that has never been adopted before for the simulation of DRT and that is particular suitable to treat the presence of corners.

The main findings of our simulations reveal that, in general, the presence of corners appearing in the particle shape is not so important in determining the DRT kinetics. 
Indeed, due to the very high corner curvature (theoretically infinite), the local mass transport coefficient $K$ is very high in comparison to that competing to the other parts of the surface. 
Thus, very rapidly, corners disappear and the particle shape tend to that of a sphere. 
This is the reason why, usually, the effect of corners on the whole dissolution phenomenon is not so relevant. 
Accordingly, our model makes it clear that the presence of small irregularities on particle surface should not matter so much in the theoretical interpretation of the DRT time course. 
On the other hand, particles shape is very important, being the particle surface/volume ($S/V$) ratio one of the most important parameters. 
Indeed, the biggest $S/V$, the faster the dissolution process occurs. In this sense, reasonably, we found a strict analogy between particle dissolution and particle melting \cite{chiarappa2017}. 
Indeed, both dissolution and melting imply network destruction with the only difference that in the first case network breakdown is due to the external solvent, while in the second it is induced by temperature. 
Thus, increasing $S/V$ -- implying the enhancement of surface molecules number -- promotes dissolution kinetics and reduces (nano) particle melting temperature. 

Similarly, particles size distribution plays a very important role in ruling dissolution since, again, it is strictly connected to the $S/V$ ratio, this time referring to all the particles and not to a specific one. 
Obviously, in real situations, the ensemble of dissolving particles is made up by particles of different shapes and different sizes. Interestingly, the developed model can account also for this real situation, so that it can safely simulate practical situations and it can be used to evaluate possible equivalence among particles ensembles characterized by particles of different shape and/or size distribution. 
Our simulations reveal that DRT kinetics not only depends on $S/V$ but it also depends on particles shape and size distribution inside each shape, i.e. the situation met in real situations.
Clearly, this is very important in the light of quality control that recurs to DRT to establish equivalence among different particles lots. 
At the same time, the model can be used to study some important surface properties such as the presence of a surface phase, differing from the bulk one, and induced by dissolution (surface re-crystallization). 
In addition, the proposed model accounts also for solid wettability providing an important generalization of the expression of the interfacial mass transfer coefficient $k_m$ \cite{abrami2022EJPB}. 

Finally, we do believe that the proposed model represents a very important theoretical tool to interpret the experimental evidences coming from the developing experimental set up aimed to visualize particle shape variation upon dissolution \cite{laitinen2010, stukelj2020, lausch2024}. 
This, clearly, will imply to move from a two-dimensional model to a three-dimensional one.



\section*{Acknowledgments}
E.C., F.L.I., and G.P.\ are members of the INdAM--GNCS research group.

\section*{Funding}
E.C.\ and F.L.I.\ was funded by INdAM--GNCS Project entitled ``Numerical modeling and high-performance computing approaches for multiscale models of complex systems'' (PI: M.\ Menci), CUP E53C23001670001.

\bibliographystyle{unsrt} 
\bibliography{biblio}

\begin{thebibliography}{10}

\bibitem{siepmann2013}
J.~Siepmann and F.~Siepmann.
\newblock Mathematical modeling of drug dissolution.
\newblock {\em International Journal of Pharmaceutics}, 453(1):12--24, 2013.

\bibitem{hasa2014}
D.~Hasa, B.~Perissutti, D.~Voinovich, M.~Abrami, R.~Farra, S.~M. Fiorentino,
  G.~Grassi, and M.~Grassi.
\newblock Drug nanocrystals: theoretical background of solubility increase and
  dissolution rate enhancement.
\newblock {\em Chemical and Biochemical Engineering Quarterly}, 28(3):247--258,
  2014.

\bibitem{parmar2018}
A.~Parmar and S.~Sharma.
\newblock Engineering design and mechanistic mathematical models: standpoint on
  cutting edge drug delivery.
\newblock {\em {TrAC} Trends in Analytical Chemistry}, 100:15--35, 2018.

\bibitem{kwan1997}
K.~C. Kwan.
\newblock Oral bioavailability and first-pass effects.
\newblock {\em Drug Metabolism and Disposition}, 25(12):1329--1336, 1997.

\bibitem{naidu2008}
R.~Naidu, K.~T. Semple, M.~Megharaj, A.~L. Juhasz, N.~S. Bolan, S.~K. Gupta,
  B.~E. Clothier, and R.~Schulin.
\newblock Bioavailability: definition, assessment and implications for risk
  assessment.
\newblock {\em Developments in Soil Science}, 32:39--51, 2008.

\bibitem{loftsson2010}
T.~Loftsson and M.~E. Brewster.
\newblock Pharmaceutical applications of cyclodextrins: basic science and
  product development.
\newblock {\em Journal of Pharmacy and Pharmacology}, 62(11):1607--1621, 2010.

\bibitem{davis2018}
M.~Davis and G.~Walker.
\newblock Recent strategies in spray drying for the enhanced bioavailability of
  poorly water-soluble drugs.
\newblock {\em Journal of Controlled Release}, 269:110--127, 2018.

\bibitem{gigliobianco2018}
M.~R. Gigliobianco, C.~Casadidio, R.~Censi, and P.~Di~Martino.
\newblock Nanocrystals of poorly soluble drugs: drug bioavailability and
  physicochemical stability.
\newblock {\em Pharmaceutics}, 10(3):134, 2018.

\bibitem{bertoni2019}
S.~Bertoni, B.~Albertini, and N.~Passerini.
\newblock Spray congealing: An emerging technology to prepare solid dispersions
  with enhanced oral bioavailability of poorly water soluble drugs.
\newblock {\em Molecules}, 24(19):3471, 2019.

\bibitem{hixson1931a}
A.~W. Hixson and J.~H. Crowell.
\newblock Dependence of reaction velocity upon surface and agitation {I}.
  {T}heoretical consideration.
\newblock {\em Industrial \& Engineering Chemistry}, 23(8):923--931, 1931.

\bibitem{hixson1931b}
A.~W. Hixson and J.~H. Crowell.
\newblock Dependence of reaction velocity upon surface and agitation {II}.
  {E}xperimental procedure in study of surface.
\newblock {\em Industrial \& Engineering Chemistry}, 23(9):1002--1009, 1931.

\bibitem{hixson1931c}
A.~W. Hixson and J.~H. Crowell.
\newblock Dependence of reaction velocity upon surface and agitation {III}.
  {E}xperimental procedure in study of agitation.
\newblock {\em Industrial \& Engineering Chemistry}, 23(10):1160--1169, 1931.

\bibitem{niebergall1963}
P.~J. Niebergall, G.~Milosovich, and J.~E. Goyan.
\newblock Dissolution rate studies {II}: Dissolution of particles under
  conditions of rapid agitation.
\newblock {\em Journal of Pharmaceutical Sciences}, 52(3):236--241, 1963.

\bibitem{pedersen1975a}
P.~V. Pedersen and K.~F. Brown.
\newblock Dissolution profile in relation to initial particle distribution.
\newblock {\em Journal of Pharmaceutical Sciences}, 64(7):1192--1195, 1975.

\bibitem{pedersen1975b}
P.~V. Pedersen and K.~F. Brown.
\newblock Size distribution effects in multiparticulate dissolution.
\newblock {\em Journal of Pharmaceutical Sciences}, 64(12):1981--1986, 1975.

\bibitem{pedersen1977}
P.~V. Pedersen and K.~F. Brown.
\newblock General class of multiparticulate dissolution models.
\newblock {\em Journal of Pharmaceutical Sciences}, 66(10):1435--1438, 1977.

\bibitem{pedersen1978}
P.~V. Pedersen and J.~W. Myrick.
\newblock Versatile kinetic approach to analysis of dissolution data.
\newblock {\em Journal of Pharmaceutical Sciences}, 67(10):1450--1455, 1978.

\bibitem{grassi2007}
M.~Grassi, G.~Grassi, R.~Lapasin, and I.~Colombo.
\newblock {\em Understanding drug release and absorption mechanisms: a physical
  and mathematical approach}.
\newblock CRC Press, Boca Raton, USA, 2007.

\bibitem{guo2018}
N.~Guo, B.~Hou, N.~Wang, Y.~Xiao, J.~Huang, Y.~Guo, S.~Zong, and H.~Hao.
\newblock In situ monitoring and modeling of the solution-mediated polymorphic
  transformation of rifampicin: from form {II} to form {I}.
\newblock {\em Journal of Pharmaceutical Sciences}, 107(1):344--352, 2018.

\bibitem{abrami2020ADMET}
M.~Abrami, L.~Grassi, R.~Di~Vittorio, D.~Hasa, B.~Perissutti, D.~Voinovich,
  G.~Grassi, I.~Colombo, and M.~Grassi.
\newblock Dissolution of an ensemble of differently shaped poly-dispersed drug
  particles undergoing solubility reduction: mathematical modelling.
\newblock {\em ADMET and DMPK}, 8(3):297--313, 2020.

\bibitem{thormann2014}
U.~Thormann, M.~De~Mieri, M.~Neuburger, S.~Verjee, P.~Altmann, M.~Hamburger,
  and G.~Imanidis.
\newblock Mechanism of chemical degradation and determination of solubility by
  kinetic modeling of the highly unstable sesquiterpene lactone nobilin in
  different media.
\newblock {\em Journal of Pharmaceutical Sciences}, 103(10):3139--3152, 2014.

\bibitem{mosharraf1995}
M.~Mosharraf and C.~Nystr{\"o}m.
\newblock The effect of particle size and shape on the surface specific
  dissolution rate of microsized practically insoluble drugs.
\newblock {\em International Journal of Pharmaceutics}, 122(1-2):35--47, 1995.

\bibitem{hirai2017}
D.~Hirai, Y.~Iwao, S.-I. Kimura, S.~Noguchi, and S.~Itai.
\newblock Mathematical model to analyze the dissolution behavior of metastable
  crystals or amorphous drug accompanied with a solid-liquid interface
  reaction.
\newblock {\em International Journal of Pharmaceutics}, 522(1-2):58--65, 2017.

\bibitem{hsu2021}
S.-Y. Hsu and C.-W. Wu.
\newblock Effects of particle shape in mass-diffusion-controlled dissolution
  process.
\newblock {\em International Communications in Heat and Mass Transfer},
  125:105321, 2021.

\bibitem{yuan2013}
Q.~Yuan, X.~Jia, and R.~A. Williams.
\newblock Validation of a multi-component digital dissolution model for
  irregular particles.
\newblock {\em Powder Technology}, 240:25--30, 2013.

\bibitem{lausch2024}
M.~Lausch, P.~Brockmann, F.~Schmitt, B.~J.~M. Etzold, and J.~Hussong.
\newblock In-situ iron oxide particle size and shape evolution during the
  dissolution in oxalic acid.
\newblock {\em Chemical Engineering Science}, 289:119864, 2024.

\bibitem{hsu1991}
J.-P. Hsu, D.-L. Lin, and M.-J. Lin.
\newblock Dissolution of solid particles in liquids: A surface layer model.
\newblock {\em Colloids and Surfaces}, 61:35--47, 1991.

\bibitem{wang1999}
J.~Wang and D.~R. Flanagan.
\newblock General solution for diffusion-controlled dissolution of spherical
  particles. 1. {T}heory.
\newblock {\em Journal of Pharmaceutical Sciences}, 88(7):731--738, 1999.

\bibitem{salish2024}
K.~Salish, C.~So, S.~H. Jeong, H.~H. Hou, and C.~Mao.
\newblock A refined thin-film model for drug dissolution considering radial
  diffusion -- simulating powder dissolution.
\newblock {\em Pharmaceutical Research}, 41:947--958, 2024.

\bibitem{darcy2011}
D.~M. D'{A}rcy and T.~Persoons.
\newblock Mechanistic modelling and mechanistic monitoring: simulation and
  shadowgraph imaging of particulate dissolution in the flow-through apparatus.
\newblock {\em Journal of Pharmaceutical Sciences}, 100(3):1102--1115, 2011.

\bibitem{darcy2019}
D.~M. D'{A}rcy and T.~Persoons.
\newblock Understanding the potential for dissolution simulation to explore the
  effects of medium viscosity on particulate dissolution.
\newblock {\em AAPS PharmSciTech}, 20:47, 2019.

\bibitem{abrami2022EJPB}
M.~Abrami, M.~Grassi, D.~Masiello, and G.~Pontrelli.
\newblock Dissolution of irregularly-shaped drug particles: mathematical
  modelling.
\newblock {\em European Journal of Pharmaceutics and Biopharmaceutics},
  177:199--210, 2022.

\bibitem{osher1988JCP}
S.~Osher and J.~A. Sethian.
\newblock Front propagating with curvature-dependent speed: algorithms based on
  {H}amilton-{J}acobi formulations.
\newblock {\em J. Comput. Phys.}, 79(1):12--49, 1988.

\bibitem{nogami1969}
H.~Nogami, T.~Nagai, and T.~Yotsuyanagi.
\newblock Dissolution phenomena of organic medicinals involving simultaneous
  phase changes.
\newblock {\em Chemical and Pharmaceutical Bulletin}, 17(3):499--509, 1969.

\bibitem{levich1962book}
V.~G. Levich.
\newblock {\em Physicochemical Hydrodynamics}.
\newblock Prentice-Hall, Englewood Cliffs, N.J., 1962.

\bibitem{osherbook}
S.~Osher and R.~Fedkiw.
\newblock {\em Level set methods and dynamic implicit surfaces}, volume 153 of
  {\em Applied Mathematical Sciences}.
\newblock Springer-Verlag, New York, 2003.

\bibitem{sethianbook}
J.~A. Sethian.
\newblock {\em Level set methods and fast marching methods: evolving interfaces
  in computational geometry, fluid mechanics, computer vision, and material
  science}.
\newblock Cambridge University Press, New York, 1999.

\bibitem{naseem2004}
A.~Naseem, C.~J. Olliff, L.~G. Martini, and A.~W. Lloyd.
\newblock Effects of plasma irradiation on the wettability and dissolution of
  compacts of griseofulvin.
\newblock {\em International Journal of Pharmaceutics}, 269(2):443--450, 2004.

\bibitem{chiarappa2017}
G.~Chiarappa, A.~Piccolo, I.~Colombo, D.~Hasa, D.~Voinovich, M.~Moneghini,
  G.~Grassi, R.~Farra, M.~Abrami, P.~Posocco, S.~Prici, and M.~Grassi.
\newblock Exploring the shape influence on melting temperature, enthalpy, and
  solubility of organic drug nanocrystals by a thermodynamic model.
\newblock {\em Crystal Growth \& Design}, 17(8):4072--4083, 2017.

\bibitem{laitinen2010}
R.~Laitinen, J.~Lahtinen, P.~Silfsten, E.~Vartiainen, P.~Jarho, and
  J.~Ketolainen.
\newblock An optical method for continuous monitoring of the dissolution rate
  of pharmaceutical powders.
\newblock {\em Journal of Pharmaceutical and Biomedical Analysis},
  52(2):181--189, 2010.

\bibitem{stukelj2020}
J.~{\v{S}}tukelj, M.~Agopov, J.~Yliruusi, C.~J. Strachan, and S.~Svanb{\"a}ck.
\newblock Image-based dissolution analysis for tracking the surface stability
  of amorphous powders.
\newblock {\em ADMET and DMPK}, 8(4):401--409, 2020.

\bibitem{adamson1997book}
A.~W. Adamson and A.~P. Gast.
\newblock {\em Physical Chemistry of Surfaces}.
\newblock Wiley, New York, 1997.

\bibitem{tolman1949}
R.~C. Tolman.
\newblock The effect of droplet size on surface tension.
\newblock {\em The Journal of Chemical Physics}, 17(3):333--337, 1949.

\bibitem{samsonov2004}
V.~M. Samsonov, N.~Y. Sdobnyakov, and A.~N. Bazulev.
\newblock Size dependence of the surface tension and the problem of {G}ibbs
  thermodynamics extension to nanosystems.
\newblock {\em Colloids and Surfaces A: Physicochemical and Engineering
  Aspects}, 239(1-3):113--117, 2004.

\bibitem{rowlinson1982book}
J.~S. Rowlinson and B.~Widom.
\newblock {\em Molecular Theory of Capillarity}.
\newblock Clarendon Press, Oxford, 1982.

\end{thebibliography}

\clearpage
\appendix

\section*{Appendix}

\section{Particle distribution}\label{app:weibull}
It is known that the size of the particles is randomly distributed following shifted Weibull distribution \cite{abrami2020ADMET}
\begin{equation}\label{weibull}
f(x)=\frac{k}{\lambda^k}(x-x_0)^{k-1}\exp{\left[-\left(\frac{x-x_0}{\lambda}\right)^k\right]}, \qquad x>x_0.
\end{equation}

The cumulative distribution is given by
\begin{equation}\label{weibullcum}
F(x)=\int_{-\infty}^x f(\zeta)d\zeta=1-\exp{\left[-\left(\frac{x-x_0}{\lambda}\right)^k\right]}, \qquad x>x_0,
\end{equation}
and the inverse of $F$ is given by
\begin{equation}\label{weibullcuminv}
Q(y)=\lambda(-\ln(1-y))^{\frac{1}{k}}+x_0, \quad y\in[0,1].
\end{equation}

Using \eqref{weibullcuminv} one can easily extract any number $n$ of random samples $R_1,\dots,R_n$ distributed accordingly to \eqref{weibull}, simply defining
$$
R_i=Q(y_i), \qquad i=1,\ldots, n,
$$
with $y_i$ \emph{uniformly} random variables in [0,1]. 

\section{Numerical discretization}\label{app:numerics}
In this Appendix we describe in full details the numerical method we used to get simulations of Section \ref{sec:simulations}.

\subsection{Numerical grid and CFL condition}
Let us introduce a structured grid $\mathcal G$. We consider a square computational domain divided in $N\times N$ square cells of side $\Delta x$. 
Let us denote as usual by $x_{ij}$ the grid node corresponding to the point $(i\Delta x,j\Delta x)$, and by $C_{ij}$ the cell with vertices $(x_{ij},x_{i+1,j},x_{i,j+1},x_{i+1,j+1})$.

Regarding time discretization, let us denote by $\Delta t$ the time step and by $t^n=n\Delta t$ the physical time at the $n$-th time step.

A crucial observation regards the CFL condition: indeed, in such a differential problems involving hyperbolic partial differential equations like \eqref{LSequation}, the ratio between time and space steps must be chosen in such a way that
\begin{equation}\label{CFL}
    \frac{\Delta t}{\Delta x}\leq \frac{1}{\max_{\mathcal G}|v_{ij}|}.
\end{equation}
This guarantees the numerical stability.

\subsection{Shape initialization (computation of $\varphi_0$)}\label{app:inizialization}
Computing the initial condition $\u_0$ for the level-set function (under constraint \eqref{phi0}) is a problem \textit{per se}. 
The difficulty is increased by the fact that the accuracy of the following evolution of $\u$ strongly depends on the accuracy of the initial condition $\u_0$.
A typical choice for $\u_0$ is the \emph{signed distance function} from the initial curve.
To begin with, let us recall the equations in polar coordinates of the curves we have considered:
\begin{center}
\begin{tabular}{lll}
 \textit{circle} & $r(\theta)=R$, & $\theta\in[0,2\pi)$, \\
 %
 %
 \textit{superellipse} & $r(\theta)=\left(\left|\ds\frac{\cos\theta}{a}\right|^n+\left|\ds\frac{\sin\theta}{b}\right|^n\right)^{-1/n}$, &  $\theta\in[0,2\pi)$, \ $n\geq 2$ (ellipse with $n=2$).
\end{tabular}

\end{center}

Once we have this, we can easily get the points of the initial curves with arbitrary accuracy. 
Then, for any grid node, we can compute the distance between the node and the curve by finding the minimal distance between the node and any point of the curve. This step is relatively expensive but it must done only once.
The last step requires to discriminate if the point is internal or external. This can be done by solving the well known point-in-polygon problem using, e.g., the ray casting algorithm or the winding number algorithm.

\subsection{Numerical scheme for the level-set equation}
Here we recall the first-order upwind numerical scheme we have used to solve \eqref{LSequation} in two dimension \cite[Sect.\ 6.4]{sethianbook}.
As usual we denote by $\u_{ij}^n$ the approximation of the level-set function $\u$ at the grid node $(i\Delta x,j\Delta x)$ and time $n\Delta t$. Similarly for $v_{ij}^n$.

The scheme reads as
$$
\u_{ij}^{n+1}=\u_{ij}^n-\Delta t\left(\max(v_{ij}^n,0)\Lambda^n_+ + \min(v_{ij}^n,0)\Lambda^n_-\right),\qquad \forall i,j, n>0
$$
with
\begin{eqnarray*}
    \Lambda^n_+:=\sqrt{
    \max(D_x^{-}\u^n_{ij},0)^2 + 
    \min(D_x^{+}\u^n_{ij},0)^2 + 
    \max(D_y^{-}\u^n_{ij},0)^2 + 
    \min(D_y^{+}\u^n_{ij},0)^2}
    \\
    \Lambda^n_-:=\sqrt{
    \max(D_x^{+}\u^n_{ij},0)^2 + 
    \min(D_x^{-}\u^n_{ij},0)^2 + 
    \max(D_y^{+}\u^n_{ij},0)^2 + 
    \min(D_y^{-}\u^n_{ij},0)^2}
\end{eqnarray*}
and
\begin{eqnarray*}
D_x^{+}\u^n_{ij}:=\frac{\u^n_{i+1,j}-\u^n_{ij}}{\Delta x}, \qquad
D_x^{-}\u^n_{ij}:=\frac{\u^n_{ij}-\u^n_{i-1,j}}{\Delta x}, \\
D_y^{+}\u^n_{ij}:=\frac{\u^n_{i,j+1}-\u^n_{ij}}{\Delta x}, \qquad
D_y^{-}\u^n_{ij}:=\frac{\u^n_{ij}-\u^n_{i,j-1}}{\Delta x}.
\end{eqnarray*}

\subsection{Computation of $\rosc$}
Let us recall from Sect.\ \ref{sec:K} that $\rosc$ is the inverse of the curvature. In 2D, for any $x$, $y$, the curvature can be written in terms of the level-set function using the formula \cite{sethianbook},
$$
\kappa=\frac{\u_{xx}\u_y^2-2\u_x\u_y\u_{xy}+\u_{yy}\u_x^2}{(\u_x^2+\u_y^2)^{3/2}}.
$$
The derivatives of $\u$ can be approximated, e.g., by centred finite differences,
$$
\u_{x} \approx  \frac{\u_{i+1,j}-\u_{i-1,j}}{2\Delta x},
\qquad
\u_{y} \approx  \frac{\u_{i,j+1}-\u_{i,j-1}}{2\Delta x}
$$
$$
\u_{xx} \approx  \frac{\u_{i+1,j}-2\u_{i,j}+\u_{i-1,j}}{\Delta x^2},
\qquad
\u_{yy} \approx  \frac{\u_{i,j+1}-2\u_{i,j}+\u_{i,j-1}}{\Delta x^2},
$$
$$
\u_{xy} \approx  \frac{\u_{i+1,j+1}-\u_{i+1,j-1}-\u_{i-1,j+1}+\u_{i-1,j-1}}{4\Delta x^2}.
$$

\subsection{Computation of perimeter and area of the particle}
One of the most technical steps in the numerical approximation of the level-set equation \eqref{mainequation_singleparticle} is the computation of the perimeter and the area of the particle, which are needed to evaluate $\req$ and then $\int_{\u=0} K ds$ by a quadrature formula, at any time step.

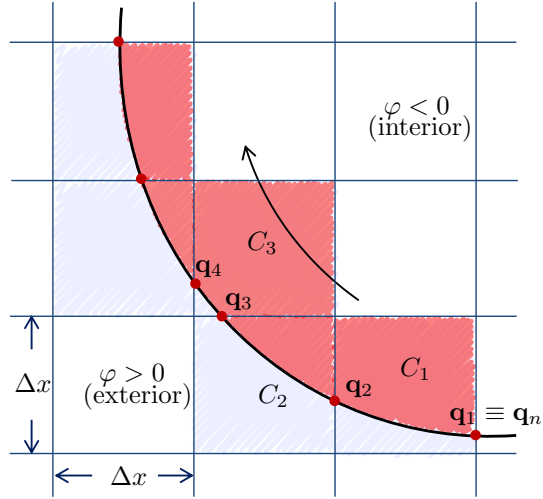
\begin{figure}[h!]
    \centering
    \begin{overpic}[width=0.6\textwidth]{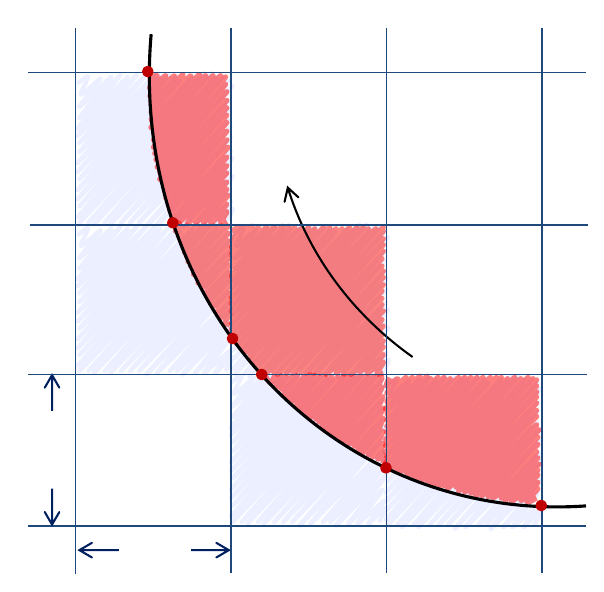}
    \put(85,17){$\Gq_1\equiv \Gq_n$} \put(66,22){$\Gq_2$} \put(44,38){$\Gq_3$} \put(38.5,44){$\Gq_4$}
    \put(76,25){$C_1$} \put(50,20){$C_2$} \put(48,48){$C_3$}
    \put(21,24){$\u>0$} \put(73,73){$\u<0$}
    \put(23,5){$\Delta x$} \put(5.5,22){$\Delta x$}
    \put(70,69){(interior)}
    \put(18,20){(exterior)}
    \end{overpic}
    \caption{Localization of the boundary of a particle on the grid.}
    \label{fig:approximated-area}
\end{figure}

The goal here is to locate the boundary of the particle by means of the level-set function $\u$. This can be done noting that the boundary divides the space in two regions, characterized by $\u>0$ (the exterior) and $\u<0$ (the interior). 
The grid cells which are crossed by the boundary are easily found since the values of $\u$ at their four vertices have different sign. By interpolation (we use linear) one can locate one point, say $\Gq_1$, of the boundary passing through a side of a grid cell, see Fig.\ \ref{fig:approximated-area}.

Once the first point is found, it is labeled as ``incoming'' in a cell (cell $C_1$ in the Figure), Then, by comparing the signs of $\u$ along the vertices of the cell, one can search for the ``outgoing'' point $\Gq_2$ (which corresponds to the incoming point in the next cell $C_2$). 
Repeating the procedure one can pass to one cell to another until the first point is reached again. 
At this point the perimeter is easily computed by summing all the lengths of the lines joining one point with the following one, i.e.
$
\sum_{i=1}^{n-1} \| \Gq_{i+1}-\Gq_{i}\|
$.

Regarding the area, the cells fully inside the particle are easily found, being characterized by having $\u<0$ at the four vertices. 
In order to increase the accuracy, one can add to the areas of the fully contained cells, the areas of the internal part of the cells crossed by the particle boundary (red areas in the Figure).

\section{Derivation of $k_m$}\label{app:derivation_km}

In order to derive the theoretical dependence of the interface mass
transport coefficient $k_m$ on the local curvature
radius $R$, it is possible taking advantage of the infinitesimal
variation of the internal energy $U$ competing
to the interface separating the solid and the solution phases \cite{adamson1997book}:
\begin{equation}\label{A01}
	dU=T dS+
	\sum_{i=1}^r\mu_i dn_i+
	\gamma dA+
	C_1d\kappa_1+
	C_2d\kappa_2,
\end{equation}
where -- limited to this Appendix -- $T$ and $S$ denote,
respectively, the interface temperature and entropy,
$\mu_i$ and $n_i$ denote, respectively, the
interface chemical potential and moles number competing to species $i^{th}$, 
$\gamma$ and $A$ denote, respectively, the
interface surface energy and area, 
while $\kappa_1$ and $\kappa_2$ represent, respectively, the first and the
second interface curvature, 
being $C_1$ and $C_2$ two constants. 

Dividing \eqref{A01} for $dA$ and remembering that only two
species can exist (1 = solid; 2 = solvent), we have:
\begin{equation}\label{A02}
	\frac{\partial U}{\partial A}=
	T\frac{\partial S}{\partial A}+
	\mu_1\frac{\partial n_1}{\partial A}+
	\mu_2\frac{\partial n_2}{\partial A}+
	\gamma+
	C_1\frac{\partial \kappa_1}{\partial A}+
	C_2\frac{\partial \kappa_2}{\partial A}.
\end{equation}

As $S$ and $n_i$ are extensive functions of
$A$ and increase linearly with $A$, we can set:
\begin{equation}\label{A03}
	\frac{\partial S}{\partial A}=k_S, \qquad 
	\frac{\partial n_1}{\partial A}=k_{n_1}, \qquad
	\frac{\partial n_2}{\partial A}=k_{n_2}.
\end{equation}

On the basis of \eqref{A03} and assuming that $n_2$ is negligible in comparison to $n_1$ (this assumption relies
on the usual hypothesis that the liquid phase does not pervade the solid
one in a solid-liquid equilibrium), \eqref{A02} becomes:
\begin{equation}\label{A04}
	\frac{\partial U}{\partial A}\approx
	T k_S+
	\mu_1 k_{n_1}+
	\gamma+
	C_1\frac{\partial \kappa_1}{\partial A}+
	C_2\frac{\partial \kappa_2}{\partial A}.
\end{equation}

If we further suppose that $\kappa_1\approx \kappa_2$ (we locally approximate the surface by a sphere), it follows:
\begin{equation}\label{A05}
	\frac{\partial \kappa_1}{\partial A}\approx
	\frac{\partial \kappa_2}{\partial A}=
	\frac{\partial\left(\frac1R\right)}{\partial(4\pi R^2)}=
	-\frac{1}{8\pi R^3}.
\end{equation}

Inserting \eqref{A05} into \eqref{A04} leads to:
\begin{equation}\label{A06}
	\frac{\partial U}{\partial A}\approx
	T k_S+
	\mu_1 k_{n_1} +
	\gamma-
	\frac{(C_1+C_2)}{8\pi R^3}.
\end{equation}

As equilibrium conditions require $dU=0$, we have:
\begin{equation}\label{A07}
	T k_S+
	\mu_1 k_{n_1} +
	\gamma-
	\frac{(C_1+C_2)}{8\pi R^3}=0.
\end{equation}

Solving \eqref{A07} for $\mu_1$ leads to:
\begin{equation}\label{A08}
	\mu_1=
	\frac{(C_1+C_2)}{8\pi k_{n_1}}\frac{1}{R^3}-
	\frac{\gamma}{k_{n_1}}-
	\frac{T k_S}{k_{n_1}}.
\end{equation}

In the light of the Tolman equation \cite{tolman1949}:
\begin{equation}\label{A09}
	\gamma=\gamma_\infty\frac{R}{R+2d_T}
\end{equation}
where $\gamma_\infty$ is the solid-liquid surface energy for a flat
surface (vanishing curvature) and $d_T$ is the Tolman length, whose order of magnitude should correspond to the
diameter of the molecules constituting the curved surface \cite{samsonov2004}, but it is usually assumed to be 1/3 of the molecules diameter \cite{rowlinson1982book}, \eqref{A08} becomes:
\begin{equation}\label{A10}
	\mu_1=
	\frac{(C_1+C_2)}{8\pi k_{n_1}}\frac{1}{R^3}-
	\frac{\gamma_\infty}{k_{n_1}}\frac{R}{R+2d_T}-
	\frac{T k_S}{k_{n_1}}.
\end{equation}

As the chemical potential expresses the tendency of solid surface
molecules to leave the interface, we can argue that
$k_m$ is proportional to $\mu_1$. 
Assuming that $k_m$ gets its minimum value (that corresponding to
a plane surface, $R\to\infty$) when $\mu_1=0$, \eqref{A10}
seems to suggest that a reasonable $k_m$ dependence
on $R$ is given by:
\begin{equation}\label{A11final}
	k_{m} = k_m^\infty \left( \frac{\alpha}{R^3} + \frac{R}{R + 2d_T} \right)
\end{equation}
where $\alpha$ is an unknown parameter to be determined.

\end{document}